\documentclass[12pt,a4paper]{article}
\usepackage{amssymb}
\usepackage{amsmath}
\usepackage{latexsym}
\usepackage{amsthm}

\newcommand{\V}{\textup{Vol}(4_1)}
\newcommand{\ve}{\varepsilon}

\newtheorem{thm}{Theorem}
\newtheorem{lem}[thm]{Lemma}
\newtheorem{cor}[thm]{Corollary}
\newtheorem{prop}[thm]{Proposition}
\theoremstyle{definition}
\newtheorem*{remark}{Remark}
\theoremstyle{definition}
\newtheorem{definition}{Definition}

\usepackage{graphicx}
\usepackage{subfigure}
\usepackage{color}

\usepackage{xpatch}
\xpatchcmd{\proof}{\itshape}{\normalfont\proofnameformat}{}{}
\newcommand{\proofnameformat}{}

\newcommand\blfootnote[1]{%
  \begingroup
  \renewcommand\thefootnote{}\footnote{#1}%
  \addtocounter{footnote}{-1}%
  \endgroup
}

\begin{document}

\renewcommand{\proofnameformat}{\bfseries}

\begin{center}
{\Large\textbf{Maximizing Sudler products via Ostrowski expansions and cotangent sums}}


\vspace{5mm}

\textbf{Christoph Aistleitner$^1$ and Bence Borda$^{1,2}$}

\vspace{5mm}

{\footnotesize $^1$Graz University of Technology

Institute of Analysis and Number Theory

Steyrergasse 30, 8010 Graz, Austria

\vspace{5mm}

$^2$Alfr\'ed R\'enyi Institute of Mathematics

Re\'altanoda utca 13--15, 1053 Budapest, Hungary

Email: \texttt{aistleitner@math.tugraz.at} and \texttt{borda@math.tugraz.at}}
\end{center}

\vspace{5mm}

\begin{abstract}
There is an extensive literature on the asymptotic order of Sudler's trigonometric product $P_N (\alpha) = \prod_{n=1}^N |2 \sin (\pi n \alpha)|$ for fixed or for ``typical'' values of $\alpha$. In the present paper we establish a structural result, which for a given $\alpha$ characterizes those $N$ for which $P_N(\alpha)$ attains particularly large values. This characterization relies on the coefficients of $N$ in its Ostrowski expansion with respect to $\alpha$, and allows us to obtain very precise estimates for $\max_{1 \le N \leq M} P_N(\alpha)$ and for $\sum_{N=1}^M P_N(\alpha)^c$ in terms of $M$, for any $c>0$. Furthermore, our arguments give a natural explanation of the fact that the value of the hyperbolic volume of the complement of the figure-eight knot appears generically in results on the asymptotic order of the Sudler product and of the Kashaev invariant.
\end{abstract}

\section{Introduction and statement of results}

During the last decades many authors have studied the asymptotic order of the so-called Sudler product
\begin{equation} \label{sudler}
P_N(\alpha) = \prod_{n=1}^N |2 \sin (\pi n \alpha)| ,
\end{equation}
either on average (with respect to $\alpha$) or for particular values of $\alpha$. It is known that the order of \eqref{sudler} for a fixed value of $\alpha$ depends sensitively on the Diophantine approximation properties of $\alpha$, and in particular on the continued fraction expansion of $\alpha$. If $\alpha \in \mathbb{Q}$, then clearly the product \eqref{sudler} vanishes for all sufficiently large $N$, so for the asymptotic analysis we can restrict ourselves to the case when $\alpha$ is irrational.

Of particular interest is the case when $\alpha$ is a quadratic irrational, which means that the continued fraction expansion of $\alpha$ is eventually periodic. A remarkable result was recently obtained by Grepstad, Kaltenb\"ock and Neum\"uller \cite{GKN}, who proved that for the golden mean $\phi =(1+\sqrt{5})/2$,
$$
0< \liminf_{N \to \infty} P_N(\phi) < \infty ,
$$
thereby solving a long-standing problem of Erd\H os and Szekeres \cite{ESZ}. In \cite{ATZ} this result was complemented by 
$$
0 < \limsup_{N \to \infty} \frac{P_N(\phi)}{N} < \infty,
$$
so the asymptotic order of $P_N(\phi)$ is completely understood. Interestingly, there is a transition in the behavior of quadratic irrationals whose continued fraction expansion is of the particularly simple form $\alpha = [0;\overline{a}]$ as the value of $a$ increases (here and in the sequel, the overline denotes period): it turns out that $\liminf_{N \to \infty} P_N(\alpha) > 0$ as long as $a \leq 5$, while $\liminf_{N \to \infty} P_N(\alpha) = 0$ when $a \geq 6$. A similar characterization applies to $\limsup_{N \to \infty} P_N(\alpha)/N < \infty$. A generalization of such a criterion to more general quadratic irrationals can be found in \cite{GNZ}.

In \cite{LU} Lubinsky proved that for any badly approximable $\alpha$,
$$
N^{-c_1} \ll_{\alpha} P_N(\alpha) \ll_{\alpha} N^{c_2}
$$
with some $c_1, c_2 \ge 0$, and asked for the smallest possible constants $c_1=c_1(\alpha)$ and $c_2=c_2(\alpha)$ for which this holds.\footnote{Throughout the paper we write $x_N \ll y_N$ or $x_N=O(y_N)$ when $|x_N| \le C y_N$ with some appropriate constant $C>0$. All implied constants are universal unless the opposite is explicitly indicated by a subscript; e.g.\ $x_N \ll_{\alpha} y_N$ and $x_N=O_{\alpha}(y_N)$ mean that the implied constant may depend on $\alpha$.} As noted in \cite{AB}, for any badly approximable $\alpha$ we have $c_2(\alpha) = c_1(\alpha) +1$. From what was said above for $\alpha = [0;\overline{a}]$, we have $c_1(\alpha)=0$ and $c_2(\alpha)=1$ for $a \in \{1,2,3,4,5\}$, but in general it seems to be very difficult to calculate the values of these two constants. In \cite{AB} it was shown that for any quadratic irrational $\alpha$,
\begin{equation} \label{c2c1_est}
c_2(\alpha) = c_1(\alpha) + 1 = \frac{\textup{Vol} (4_1)}{4 \pi} \cdot \frac{a_{\textup{avg}}(\alpha)}{\log \lambda(\alpha)} + O \left( \frac{1+\log A(\alpha)}{\log \lambda(\alpha)} \right). 
\end{equation}
Here 
\begin{equation} \label{vol41}
\textup{Vol}(4_1) = 4 \pi \int_0^{5/6} \log (2 \sin (\pi x)) \, \mathrm{d}x \approx 2.02988
\end{equation}
is the hyperbolic volume of the complement of the figure-eight knot (more on this below; here and in the sequel, ``$4_1$'' is the Alexander--Briggs notation for the figure-eight knot), $a_{\textup{avg}}(\alpha)= \lim_{k \to \infty} (a_1 + \dots + a_k)/k$ denotes the average of the partial quotients within a period, $\lambda(\alpha)$ is the (easily computable) number for which the convergents $p_k/q_k = [a_0;a_1,a_2, \dots, a_k]$ satisfy $\log q_k \sim (\log \lambda(\alpha)) k$ as $k \to \infty$, and $A(\alpha )=\max_{k \ge 1} a_k$ is the maximum of the partial quotients. The key purpose of the present paper is to obtain a significantly improved version of \eqref{c2c1_est}, and to give a structural description of those values of $N$ for which $P_N(\alpha)$ attains particularly large resp.\ small values. These two aims are very closely related; roughly speaking, knowing the particular structure of those $N$ which lead to extreme values of $P_N(\alpha)$ allows us to obtain improved estimates on $c_1 (\alpha )$ and $c_2 (\alpha )$, since the structural information allows a refined analysis of the terms that control $c_1 (\alpha )$ and $c_2 (\alpha )$. The ``structure'' of $N$ which is alluded to here is a particular structure of the coefficients in its Ostrowski representation, which is a numeration system for integers based on the continued fraction denominators of $\alpha$. Very roughly speaking, it turns out that $P_N(\alpha)$ is particularly large resp.\ small if the Ostrowski coefficients of $N$ are all $5/6$ resp.\ $1/6$ of their maximal possible size; this fact will also give a natural explanation for the appearance of the constant $\textup{Vol}(4_1)$ in formula \eqref{c2c1_est} above, as well as in many other related formulas such as those in \cite{BD}.

Before presenting our results, we note some connections to other areas of mathematics. Early investigations of the Sudler product were carried out by Erd\H os and Szekeres \cite{ESZ} and Sudler \cite{sudler} around 1960. Since then such products have appeared in various contexts, including partition functions, KAM theory, $q$-series, Pad\'e approximations, and the analytic continuation of Dirichlet series. In particular, pointwise upper bounds for Sudler products at quadratic irrationals played a crucial role in Lubinsky's \cite{LU2} counter\-example to the Baker--Gammel--Wills conjecture, where they were used to bound the Taylor series coefficients of the Rogers--Ramanujan function. Pointwise upper bounds for Sudler products also played a key role in the solution of the ``Ten Martini Problem'' by Avila and Jitomirskaya \cite{AJ}.  However, rather than giving an exhaustive list of appearances of such products we refer the reader to \cite{ALPST,KT,VM} and the references therein. Note that the Sudler product can be written using the $q$-Pochhammer symbol as 
$$
P_N(\alpha) = |(q;q)_N| = |(1 -q)(1-q^2) \cdots (1-q^N)| \quad \text{with} \quad q = e^{2 \pi i \alpha}.
$$
Compare this with the definition of the so-called Kashaev invariant of the figure-eight knot, given by
\begin{equation} \label{kash_def}
\mathcal{J}_{4_1,0} (q) = \sum_{N=0}^\infty |(q;q)_N|^2.
\end{equation} 
This series is convergent if and only if $\alpha$ is rational (i.e.\ $q$ is a root of unity). The Kashaev invariant is a quantum knot invariant arising from the colored Jones polynomial, and the figure-eight knot is the simplest hyperbolic knot. For more background we refer to \cite{BD} and the references given there. Here we only note that the Kashaev invariant can be written as a sum of squares of Sudler products; however, by its very nature the Kashaev invariant is only interesting when $\alpha$ is rational (since otherwise the series diverges), while the asymptotic order of the Sudler product is only interesting when $\alpha$ is irrational (since otherwise the product vanishes for all sufficiently large indices). However, it is possible to approximate the value of the Sudler product at an irrational $\alpha$ by the value at a rational number close to $\alpha$ (such as a continued fraction approximation to $\alpha$), thereby switching from Kashaev invariants to Sudler products and vice versa; see \cite{AB} for a precise statement. Generally, this connection suggests to study the asymptotic order of expressions of the form
$$
\left(\sum_{N=1}^M P_N(\alpha)^2 \right)^{1/2}, \qquad \text{or, more generally}, \qquad \left(\sum_{N=1}^M P_N(\alpha)^c\right)^{1/c} \quad \text{for some $c>0$},
$$
which can be seen as describing the average order of $P_N(\alpha)$ with respect to $N$. With this notation the problem concerning the upper asymptotic order of $P_N (\alpha )$ corresponds to the maximum norm, that is, to the case $c=\infty$. In connection with a problem posed by Bettin and Drappeau \cite{BD}, in \cite{AB} we settled the case when $\alpha$ is a  quadratic irrational, showing that in this case for any real $c>0$,
\begin{equation} \label{K_c_def}
\log \left( \sum_{N=0}^{q_k-1} P_N(\alpha)^c \right)^{1/c} \sim K_c(\alpha ) k \quad \textrm{as } k \to \infty
\end{equation}
and
$$
\log \max_{0 \leq N < q_k} P_N(\alpha) \sim K_{\infty} (\alpha ) k \quad \textrm{as } k \to \infty
$$
with some constants $K_c(\alpha)$, $K_{\infty} (\alpha ) >0$. We repeat that $K_2(\alpha)$ in \eqref{K_c_def} is closely related to the Kashaev invariant as defined in \eqref{kash_def}. Note also that the constants in the question of Lubinsky can be expressed as $c_2(\alpha ) = c_1 (\alpha ) +1 = K_{\infty}(\alpha )/\log \lambda (\alpha )$, where $\lambda (\alpha )>1$ is the same easily computable constant as in \eqref{c2c1_est}. However, in \cite{AB} it remained open whether $K_c(\alpha)$ actually depends on $c$ or not. This is related to the question whether $P_N(\alpha)$ is exceptionally large only for a very small number of indices $N$, so that the sum in \eqref{K_c_def} is essentially dominated by a small number of summands which are of extremal size. In this paper we prove that this is not the case, and that (under certain technical assumptions) the overall order of the sum $\sum_{N=0}^{q_k-1} P_N(\alpha)^c$ is not caused by a small number of exceptionally large summands. In particular, $K_c(\alpha )$ does indeed depend on $c$.

We close this discussion by noting that the Kashaev invariant features prominently in Zagier's \cite{ZA} seminal paper on quantum modular forms, where it is introduced as being ``the most mysterious and in many ways the most interesting'' example. Zagier records certain modularity properties of the function $\mathcal{J}_{4_1,0}$, and suggests that the function $h(\alpha) = \log (\mathcal{J}_{4_1,0} (e^{2 \pi i \alpha}) / \mathcal{J}_{4_1,0} (e^{2 \pi i /\alpha}))$, relating the value of the Kashaev invariant at $\alpha$ to its value at $1 / \alpha$, appears to be continuous at irrationals. This continuity hypothesis has been driving much of the recent research in this area, but as a whole it is still widely open. See \cite{AB,BDlim,BD}.

We now state our main results. For the rest of the paper, we fix an irrational $\alpha=[a_0; a_1, a_2, \dots ]$ with convergents $p_k/q_k=[a_0;a_1, \dots, a_k]$.
\begin{thm}\label{maxpntheorem} Assume that
\begin{equation}\label{logak/ak+1}
\frac{\log a_k}{a_{k+1}} \le T \quad \textrm{for all } k \ge k_0
\end{equation}
with some constants $k_0, T \ge 1$. Let $N=\sum_{k=0}^{K-1} b_k q_k$ be the Ostrowski expansion of a non-negative integer, and set 
\begin{equation} \label{n*def}
N^*=\sum_{k=0}^{K-1} b_k^* q_k, \qquad \text{with} \qquad b_k^*=\lfloor (5/6) a_{k+1} \rfloor. 
\end{equation}
Then
\[ \log P_N(\alpha ) = \log P_{N^*}(\alpha ) - \sum_{k=0}^{K-1} d_k(N) + O_T \left( \sum_{k=1}^K \frac{1}{a_k} \right) + O_{\alpha}(1) \]
with some $d_k(N)$ satisfying the following for all $0 \le k \le K-1$:
\begin{enumerate}
\item[(i)] $d_k(N) \ge 0.2326 (b_k- b_k^*)^2/a_{k+1}$ with equality if $b_k=b_k^*$;
\item[(ii)] if $b_k \le 0.99 a_{k+1}$, then
\begin{equation}\label{dkNformula}
d_k (N) = a_{k+1} \int_{b_k/a_{k+1}}^{b_k^*/a_{k+1}} \log |2 \sin (\pi x)| \, \mathrm{d}x + O_T \left( \frac{|b_k-b_k^*|}{a_{k+1}} + I_{\{ b_k \le 0.01 a_{k+1} \}} \log a_{k+1} \right) .
\end{equation}
\end{enumerate}
\end{thm}

\begin{remark} In the formula above $I_{\{ b_k \le 0.01  a_{k+1}\}}$ denotes the indicator of $b_k \le 0.01 a_{k+1}$. Formula \eqref{dkNformula} gives the precise asymptotics of $d_k(N)$ in the regime $b_k-b_k^* \approx a_{k+1}$. Using a first order Taylor approximation of $\log |2 \sin (\pi x)|$ around $x=5/6$, we immediately deduce from \eqref{dkNformula} that
\begin{equation} \label{dkequ}
d_k(N) = \frac{\pi \sqrt{3}}{2} \cdot \frac{(b_k-b_k^*)^2}{a_{k+1}} + O_T \left( \frac{|b_k-b_k^*|}{a_{k+1}} + \frac{|b_k-b_k^* |^3}{a_{k+1}^2} \right) ,
\end{equation}
yielding the precise asymptotics in the regime $b_k-b_k^*=o(a_{k+1})$.
\end{remark}

Theorem \ref{maxpntheorem} asserts that $P_N(\alpha)$ is particularly large when $N = N^*$, and that an integer $N$ whose Ostrowski expansion deviates significantly from that of $N^*$ will lead to much smaller values of $P_N(\alpha)$. The magnitude of $P_N(\alpha )/P_{N^*}(\alpha )$ is quantified in terms of the ``distance'' between the Ostrowski expansions of $N$ and $N^*$. As simple illustrative examples we mention that Theorem \ref{maxpntheorem} with $T=1$ applies to $\alpha =[0;\overline{a}]$, and also to well approximable irrationals with $a_1 \le a_2 \le \cdots \le a_k \to \infty$. Note that we do not claim that the maximum is attained at precisely $N^*$; however, for example for $\alpha =[0;\overline{a}]$ it follows that the Ostrowski coefficients of the integer at which the maximum $\max_{0 \le N < q_K} P_N (\alpha )$ is attained satisfy $b_k=(5/6)a+O(1/\varepsilon )$ for all but $\le \varepsilon K$ indices $0 \le k \le K-1$.

The significance of the value $5/6$ in our definition of $N^*$ in \eqref{n*def} is that it is a solution of the equation $|2\sin (\pi x)| = 1$. From the proofs it will become visible that choosing a value of $b_k$ smaller than $(5/6)a_{k+1}$ essentially means missing out on potential factors which exceed 1, while choosing $b_k$ larger than $(5/6)a_{k+1}$ essentially leads to extra factors which are smaller than 1; clearly both effects are counterproductive if our aim is to maximize $P_N (\alpha)$. The heuristic reasoning underpinning all the constructions and results in the present paper will be described in some detail in Section \ref{sec_heuristic} below. The value $0.2326$ in (i) is explained by
\[ \frac{1}{(5/6)^2} \int_0^{5/6} \log |2 \sin (\pi x)| \, \mathrm{d}x = \frac{9 \V}{25 \pi} = 0.23260748\dots . \]
Any constant less than $9 \V /(25 \pi)$ would work; the sharpness of this value is easily seen by letting $b_k/a_{k+1} \to 0$ in \eqref{dkNformula}. The values $0.99$ resp.\ $0.01$ in (ii), on the other hand, are basically accidental; any constants $C<1$ resp.\ $C>0$ would work, with the implied constants depending also on the choice of $C$. The reason why we have to stay away from $x=0$ and $x=1$ is that the function $\log |2 \sin (\pi x)|$ has singularities there.

Condition \eqref{logak/ak+1} is related to the behavior of a cotangent sum, see Section \ref{cotangentsection}. Probably  this condition could be relaxed in some way, but it seems very difficult to obtain a version of Theorem \ref{maxpntheorem} without any regularity assumption on the relative size of the partial quotients, since for a number $\alpha$ whose partial quotients are of very different orders of magnitude  the ``optimal'' Ostrowski coefficients $b_k^*$ should depend on $a_1, a_2, \dots, a_{k+1}$ in a more complicated way than the one suggested by \eqref{n*def}; cf.\ also Figure \ref{fig:extreme} below.

Formulas \eqref{dkNformula} and \eqref{dkequ} allow us to give precise estimates for the number of integers $0 \le N < q_K$ for which $P_N(\alpha )$ is particularly large. This is stated in Theorem \ref{Lpnormtheorem} below. The value $0.01$ in the statement of the theorem could of course again be replaced by any $C>0$, with the implied constants depending also on $C>0$.
\begin{thm}\label{Lpnormtheorem} Assume that \eqref{logak/ak+1} holds and let $N^*$ be defined as in \eqref{n*def}. Then for any real $c \ge 0.01$,
\begin{equation}\label{Lpnorm}
\begin{split} & \log \left( \sum_{N=0}^{q_K-1} P_N(\alpha )^c \right)^{1/c} \\ = &\log P_{N^*}(\alpha) + \frac{1}{2c} \sum_{k=1}^K \log \frac{2a_k}{\sqrt{3}c} \\&+ O_T \left( \sum_{k=1}^K \left( \frac{\log^{1/2} (a_k/c+2)}{c^{1/2} a_k^{1/2}} +\frac{\log^{3/2} (a_k/c+2)}{c^{3/2} a_k^{1/2}} +\frac{1}{a_k} \right) \right) + O_{\alpha} (1) , \end{split}
\end{equation}
and
\begin{equation}\label{maxnorm}
\log \max_{0 \le N < q_K} P_N(\alpha ) = \log P_{N^*} (\alpha ) + O_T \left( \sum_{k=1}^K \frac{1}{a_k} \right) + O_{\alpha}(1) .
\end{equation}
\end{thm}

Our third result shows that because of the particular structure of $N^*$, we can calculate the value of $P_{N^*}(\alpha)$ up to a very high precision. 
\begin{thm}\label{pn*theorem} Assume that \eqref{logak/ak+1} holds, and let $N^*$ be defined as in \eqref{n*def}. Then
\[ \log P_{N^*} (\alpha ) = \frac{\V}{4 \pi} \sum_{k=1}^{K} a_k + \frac{1}{2} \sum_{k=1}^{K} \log a_k + O_T \left( \sum_{k=1}^{K} \frac{1 +\log (a_k a_{k+1})}{a_{k+1}} \right) +O_{\alpha} (1) . \]
\end{thm}

Let us now compare the results obtained here with the previously known best results. Consider first $\alpha =[0;\overline{a}]$. In \cite{AB} we proved that for any $0<c\le \infty$, the constant $K_c(\alpha)$ defined in \eqref{K_c_def} satisfies
\[ K_c (\alpha ) = \frac{\V}{4 \pi} a +O \left( \max \{ 1,1/c \} (1+\log a) \right) , \]
with the dependence on $c$ hidden in the error term. Taking the asymptotics as $K \to \infty$ in Theorems \ref{Lpnormtheorem} and \ref{pn*theorem}, we immediately obtain the improvement
\[ \begin{split} K_c (\alpha ) = &\frac{\V}{4 \pi} a + \frac{1}{2} \log a + \frac{1}{2c} \log \frac{2a}{\sqrt{3}c} \\ &+ O \left( \frac{\log^{1/2} (a/c+2)}{c^{1/2} a^{1/2}} +\frac{\log^{3/2} (a/c+2)}{c^{3/2} a^{1/2}} +\frac{1+\log a}{a} \right) , \qquad 0.01 \le c \le \infty . \end{split} \]
Note that the dependence on $c$ is visible in the regime $c \ll \frac{a \log \log (a+2)}{\log (a+2)}$; above this threshold the term $\frac{1}{2c} \log \frac{2a}{\sqrt{3}c}$ is negligible compared to the error term $(1+\log a)/a$, and $K_c (\alpha )$ becomes indistinguishable from
\[ K_{\infty}(\alpha ) = \frac{\V}{4 \pi} a + \frac{1}{2} \log a + O \left( \frac{1+\log a}{a} \right) . \]
As for the question of Lubinsky, the previously known best result $c_2(\alpha )=c_1(\alpha )+1=\frac{\V}{4 \pi} \cdot \frac{a}{\log a}+O(1)$ from \cite{AB} is improved to
\[ c_2 (\alpha ) = c_1 (\alpha ) +1 = \frac{\V}{4 \pi} \cdot \frac{a}{\log a} + \frac{1}{2} + O \left( \frac{1}{a} \right) , \qquad a \ge 2. \]
Theorems \ref{Lpnormtheorem} and \ref{pn*theorem} give similar improvements for more general badly approximable irrationals whose partial quotients are roughly of the same order of magnitude; this is measured by the parameter $T \ge 1$ in \eqref{logak/ak+1}.

We also obtain improvements for certain well approximable irrationals. It is known \cite{AB,BD} that if the average partial quotient $(a_1+\cdots +a_k)/k \to \infty$, then under some mild additional assumptions on $\alpha$ for any real $c>0$ we have
\[ \log \left( \sum_{N=0}^{q_k-1} P_N(\alpha )^c \right)^{1/c} \sim \frac{\V}{4 \pi} (a_1+\cdots +a_k) \qquad \textrm{as } k \to \infty , \]
and
\[ \log \max_{0 \le N < q_k} P_N (\alpha ) \sim \frac{\V}{4 \pi} (a_1+\cdots +a_k) \qquad \textrm{as } k \to \infty . \]
Theorems \ref{Lpnormtheorem} and \ref{pn*theorem} improve these under condition \eqref{logak/ak+1} by identifying logarithmic correction terms.

Finally, we mention that maximizing and minimizing $P_N (\alpha )$ are, in a sense, equivalent problems. Indeed, we observed in \cite{AB} that for an arbitrary irrational $\alpha$ and any $0 \le N < q_K$ we have
\[ \log P_N (\alpha ) + \log P_{q_K-N-1} (\alpha ) = \log q_K +O \left( \frac{1+\log \max_{1 \le k \le K}a_k}{a_{K+1}} \right) . \]
Hence $P_N (\alpha )$ is particularly small when $P_{q_K-N-1} (\alpha )$ is particularly large, and vice versa. In particular, $P_N (\alpha )$ is particularly small when $N= \sum_{k=0}^{K-1} \lfloor (1/6) a_{k+1} \rfloor q_k$, and Theorems \ref{maxpntheorem}, \ref{Lpnormtheorem} and \ref{pn*theorem} have straightforward analogues with maximum replaced by minimum.

Before coming to the more technical parts, we briefly lay out the further content of this paper. In Section \ref{shiftedsection} we introduce a perturbed version of the Sudler product, which allows a decomposition of a full product $P_N(\alpha)$ into sub-products all of which have a number of factors which is a continued fraction denominator of $\alpha$, thereby naturally bringing into play the Ostrowski expansion of $N$. In Section \ref{sec_heuristic} we give a detailed heuristic sketch of how this decomposition leads to Theorems \ref{maxpntheorem}, \ref{Lpnormtheorem}, \ref{pn*theorem}. In particular it will become clear how the constant $5/6$ in the definition of $N^*$ and how the constant $\V$ in the conclusion of the theorems arise. A key ingredient (in the heuristic as well as in the actual proofs) is the fact that the shifted products $P_{q_k}$ have a limiting behavior, in an appropriate sense. This has been experimentally observed in \cite{ATZ}, and in the present paper we give proofs for this fact which is stated as Theorems \ref{quadraticlimitfunction} and \ref{wellapproximablelimitfunction}  in Section \ref{sec_limitfunction}. Section \ref{sec_approx_shift} contains approximation formulas for shifted Sudler products, and in particular Proposition \ref{pqkapproximation}, which plays a central role in the proofs of the theorems. To obtain our approximation formula we introduce a certain cotangent sum, which controls an important part of the behavior of the shifted Sudler product. Such cotangent sums have a rich arithmetic structure, and we make crucial use of a reciprocity formula of Bettin and Conrey \cite{BC2}. Sections \ref{sec_proofs_1}--\ref{sec_proofs_3} contain the proofs of Theorems \ref{maxpntheorem}--\ref{pn*theorem}, respectively, and finally Section \ref{sec_proofs_4} contains the proofs of Theorems \ref{quadraticlimitfunction} and \ref{wellapproximablelimitfunction}.

\section{Shifted Sudler products}\label{shiftedsection}

Let
\[ P_N (\alpha , x) := \prod_{n=1}^N |2 \sin \left( \pi (n \alpha +x)  \right) |, \qquad \alpha, x \in \mathbb{R} \]
denote a shifted form of the Sudler product. Given a non-negative integer with Ostrowski expansion $N=\sum_{k=0}^{K-1} b_k q_k$, let us also introduce the notation
\begin{equation}\label{varepsilonkdef}
\varepsilon_k (N) := q_k \sum_{\ell =k+1}^{K-1} (-1)^{k+\ell} b_{\ell} \| q_{\ell} \alpha \|.
\end{equation}
It is then easy to see that
\begin{equation}\label{pnproductform}
P_N(\alpha ) = \prod_{k=0}^{K-1} \prod_{b=0}^{b_k-1} P_{q_k} (\alpha, (-1)^k (b q_k \| q_k \alpha \| + \varepsilon_k (N)) /q_k) ,
\end{equation}
which will serve as a fundamental tool in the proof of our results. This product form of $P_N(\alpha)$ was first used by Grepstad, Kaltenb\"ock and Neum\"uller \cite{GKN}, and later also in \cite{ATZ, GKN2, GNZ}; for a detailed proof of equation \eqref{pnproductform} see \cite[Lemma 2]{AB}. As we will see, here $-1<b q_k \| q_k \alpha \| + \varepsilon_k (N)<1$, therefore understanding the behavior of the function $P_{q_k}(\alpha , (-1)^k x/q_k)$ on the interval $(-1,1)$ will play a crucial role.

\subsection{The heuristic picture} \label{sec_heuristic}

Before we give the details of how to estimate the components of the product in \eqref{pnproductform}, we present a heuristic picture of how the factors in this product formula behave, what the significance of the Ostrowski coefficients of $N$ is, why the Sudler product is essentially maximized at numbers $N$ having all the Ostrowski coefficients at $5/6$ of their maximal possible size, and how the hyperbolic volume of the complement of the figure-eight knot as defined in \eqref{vol41} appears. Assume that $0 \le N < q_K$, so that $N$ has Ostrowski expansion $N = \sum_{k=0}^{K-1} b_k q_k$. Recall that $0 \leq b_k \leq a_{k+1}$. Very roughly, we have $q_k \|q_k \alpha \| \approx 1/a_{k+1}$. It turns out that $P_{q_k} (\alpha,(-1)^k x/q_k) \approx |2 \sin (\pi x)|$. This observation is formalized in a precise form in Proposition \ref{pqkapproximation} below; see also Figure \ref{fig:limf} and Theorems \ref{quadraticlimitfunction} and \ref{wellapproximablelimitfunction}. Thus ignoring the numbers $\ve_k (N)$ in \eqref{pnproductform} for the moment, we have
$$
P_{q_k} (\alpha, (-1)^k (b q_k \| q_k \alpha \| + \varepsilon_k (N)) /q_k) \approx |2 \sin (\pi b / a_{k+1} )|, 
$$
and so, offhandedly discarding the factor corresponding to $b=0$, we have
\begin{equation} \label{prod_here}
\prod_{b=0}^{b_k-1} P_{q_k} (\alpha, (-1)^k (b q_k \| q_k \alpha \| + \varepsilon_k (N)) /q_k) \approx \prod_{b=1}^{b_k-1} |2 \sin (\pi b / a_{k+1} )|.
\end{equation}

\begin{figure}[ht!]\label{figure1}
\begin{center}
  \includegraphics[width=0.6 \linewidth]{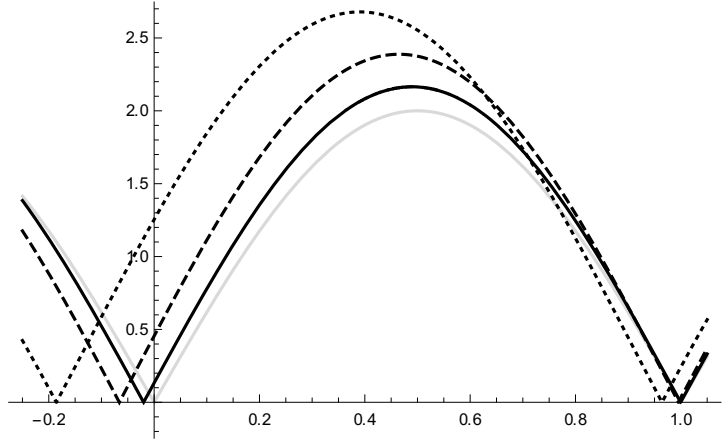}
\end{center}
  \caption{The function $P_{q_k} (\alpha,(-1)^k x/q_k)$ for $k=4$ and $\alpha = [0;\overline{a}]$, with $a=5$ (dotted), $a=15$ (dashed) and $a=50$ (solid line). The picture remains virtually identical for a larger choice of $k$. Note how the functions in the plot approach $|2 \sin (\pi x)|$ (light gray) as the value of $a$ increases. Details are given in Section \ref{sec_limitfunction} below.} 
  \label{fig:limf}
\end{figure}

Note that $b/a_{k+1} \in [0,1]$. We have $2 \sin(\pi x) \geq 1$ for $x \in [1/6,5/6]$, and $2 \sin(\pi x) \leq 1$ for $x \in [0,1/6] \cup [5/6,1]$. This suggests that in order to maximize the product in \eqref{prod_here}, we should choose $b_k \approx (5/6) a_{k+1}$, since by doing so we catch as many factors exceeding 1 while avoiding unnecessary factors smaller than 1; in other words, $P_N (\alpha )$ is essentially maximized when $N = N^*$. This heuristic also gives us a rough general approximation for the value of $P_N (\alpha)$. Using \eqref{pnproductform} and assuming that all $a_k$'s are ``large'', we roughly have
\begin{eqnarray*}
P_N (\alpha) & \approx & \prod_{k=0}^{K-1} \exp \left( \sum_{b=1}^{b_k-1} \log (2 \sin (\pi b / a_{k+1} )) \right) \\
& \approx & \exp \left( \sum_{k=0}^{K-1} a_{k+1} \int_{0}^{b_k / a_{k+1}} \log ( 2 \sin(\pi x)) \, \mathrm{d}x \right) .
\end{eqnarray*}
In particular, for $N=N^*$ when $b_k/a_{k+1} \approx 5/6$, the hyperbolic volume of the complement of the figure-eight knot naturally appears, and we have
$$
P_{N^*} ( \alpha ) \approx \exp \left( \sum_{k=0}^{K-1} a_{k+1} \int_{0}^{5/6} \log ( 2 \sin(\pi x)) \, \mathrm{d}x \right) = \exp \left( \frac{\V}{4 \pi} \sum_{k=1}^{K} a_k \right);
$$
recall the definition of $\V$ in equation \eqref{vol41}. If we want to minimize $P_N (\alpha )$ instead, the same reasoning suggests that we should choose $b_k \approx (1/6) a_{k+1}$ to catch as many factors smaller than $1$ as possible. While this heuristic serves as a good basic illustration of the behavior of the Sudler product, the actual situation clearly is much more delicate; in particular, the function $\log ( 2 \sin(\pi x))$ has singularities at $x=0$ and $x=1$, which carefully have to be taken care of.

Now let us come back to the influence of the numbers $\ve_k(N)$. As sketched above, the term $b q_k \|q_k \alpha\|$ in \eqref{prod_here} is of order roughly $b/a_{k+1}$. By \eqref{varepsilonkdef} we roughly have $|\ve_k (N)| \leq 1/a_{k+1}$, so typically the $\varepsilon_k (N)$'s are small in comparison with $b q_k \|q_k \alpha\|$. We also see in the definition given in \eqref{varepsilonkdef} that the number $\ve_k (N)$ depends on the Ostrowski coefficients $b_{k+1}, b_{k+2}, \dots$. It turns out that we cannot simply ignore the influence of the $\ve_k (N)$'s; quite on the contrary, controlling the influence of these numbers has been a key ingredient in recent work such as \cite{ATZ,GKN}, and they also play a crucial role in the present paper. In particular, the influence of the $\ve_k (N)$'s is crucial for all those factors in $P_{q_k} (\alpha, (-1)^k (b q_k \| q_k \alpha \| + \varepsilon_k (N)) /q_k) $ for which $b$ is such that $b/a_{k+1}$ is either very close to 0 or very close to 1. The punchline is the following. If a number $N$ has an Ostrowski representation which is very different from the one of $N^*$, then by the coarse argument sketched above we know that $P_N (\alpha )$ is much smaller than $P_{N^*} (\alpha )$. On the other hand, if $N$ has an Ostrowski representation which is very similar to that of $N^*$ (or in particular if $N = N^*$), then we know what the values of the $\ve_k (N)$'s are, since they depend on the Ostrowski coefficients of $N$. In other words, once we have established a structural result which controls the Ostrowski expansion of those $N$ for which  $P_N (\alpha )$ is large (Theorem \ref{maxpntheorem}), we can obtain a very precise result on the maximal asymptotic order of $P_N (\alpha )$ (combining equation \eqref{maxnorm} of Theorem \ref{Lpnormtheorem} and Theorem \ref{pn*theorem}), since control of the Ostrowski coefficients of $N$ allows us to control the numbers $\ve_k (N)$, which in turn gives us exact control of the order of $P_N (\alpha )$.

\begin{figure}[ht!]
\begin{center}
  \includegraphics[width=0.6 \linewidth]{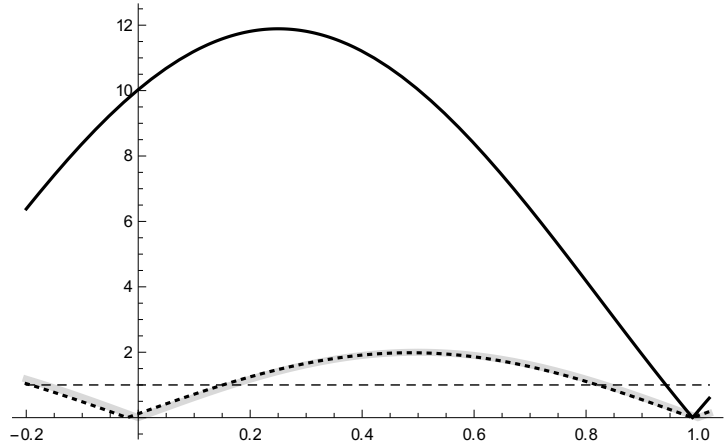}
\end{center}
  \caption{The function $P_{q_4} (\alpha,x/q_4)$ for $\alpha = [0;\overline{2,50}]$ (solid line). Note that this function is much larger than $|2 \sin (\pi x)|$ (light gray), which reflects the fact that $a_4 = 50$ is much larger than $a_3=2$. In contrast, $P_{q_5} (\alpha,-x/q_5)$ (dotted line) is virtually indistinguishable from $|2 \sin (\pi x)|$. Note that $P_{q_5} (\alpha,-x/q_5)$ crosses the line at height 1 (dashed line) near $x = 5/6$, as the initial heuristics suggested, but $P_{q_4} (\alpha,x/q_4)$ crosses this line at a much larger value of $x$ near $x = 0.95$, misleadingly suggesting a larger choice of the corresponding Ostrowski coefficient in order to maximize the Sudler product.} 
  \label{fig:extreme}
\end{figure}

There is a further important effect, which is particularly strong when $\alpha$ has some partial quotients which are very much larger than others. As Lemma \ref{Vkxlemma} and Proposition \ref{pqkapproximation} below will show, a more precise approximation for $P_{q_k}$ is
$$
P_{q_k} (\alpha, (-1)^k x/q_k) \approx |2 \sin (\pi x)| e^{(\log a_k)/a_{k+1}},
$$
where the exponential factor comes from a cotangent sum; see also Figure \ref{fig:extreme} and Sections \ref{sec_limitfunction} and \ref{cotangentsection}. If $a_k$ and $a_{k+1}$ are of similar size, then the factor $e^{(\log a_k)/a_{k+1}}$ is negligible. However, if $a_k$ is much larger than $a_{k+1}$, then this factor plays a significant role.

The heuristic sketched above suggests that in such a case the corresponding Ostrowski coefficients should be chosen significantly larger than $(5/6) a_{k+1}$, since there is a wider range of values of $x$ for which $P_{q_k} (\alpha, (-1)^k x/q_k)$ exceeds 1. However, very remarkably, this line of reasoning turns out to be wrong, and the Ostrowski coefficient maximizing the Sudler product remains at $(5/6) a_{k+1}$. The reason is that while a larger choice of $b_k$ leads to a larger value of the $k$-th factor of the Sudler product in \eqref{pnproductform}, a larger choice of $b_k$ also leads to a larger negative value of $\varepsilon_{k-1}$ which in turn leads to a smaller value of the $(k-1)$-st factor. When trying to choose a larger value of $b_k$ for some $k$ for which $(\log a_k) / a_{k+1}$ is large, then astonishingly the magnifying effect that this has on the $k$-th factor in \eqref{pnproductform} is \emph{exactly} canceled out by the corresponding de-magnifying effect on the $(k-1)$-st factor, so that overall it turns out to be better to stick with $b_k \approx (5/6) a_{k+1}$. This is a very surprising effect, which is mentioned as a ``remarkable cancellation'' in the proof of Proposition \ref{regoptprojections} (ii). We note in passing that there is a second unexpected cancellation in this paper, when the additive constant in the conclusion of Theorem \ref{pn*theorem} turns out to be zero in formula \eqref{remark_canc}. In both cases, we cannot give a convincing heuristic explanation of why these cancellations occur.

\subsection{Limit functions of shifted Sudler products} \label{sec_limitfunction}

Aistleitner, Technau and Zafeiropoulos \cite{ATZ} proved that for $\alpha=[0;\overline{a}]$ the function $P_{q_k} (\alpha, (-1)^k x/q_k)$ converges pointwise on $\mathbb{R}$ as $k \to \infty$, and gave an explicit formula for the limit function $G_{\alpha}(x)$ in the form of an infinite product. They also observed experimentally that as the value of $a$ increases the graph of $G_{\alpha}(x)$ starts to resemble that of $|2 \sin (\pi x)|$. The speed of convergence of $P_{q_k} (\alpha, (-1)^k x/q_k) \to G_{\alpha}(x)$ as $k \to \infty$ is very fast, so the graphs depicted in Figure \ref{fig:limf} for $k=4$ are practically indistinguishable from those of the corresponding limit functions $G_{\alpha} (x)$. In the present paper we develop a general framework to estimate $P_{q_k} (\alpha, (-1)^k x/q_k)$ in terms of a contangent sum; see Proposition \ref{pqkapproximation}. This in particular allows us to quantify the deviation of $P_{q_k} (\alpha, (-1)^k x/q_k)$ from $|2\sin (\pi x)|$. For the particular case of $\alpha=[0;\overline{a}]$, when passing to the limit functions $G_\alpha(x)$ by letting $k \to \infty$, we obtain
\begin{equation}\label{Galphalimitfunction}
\begin{split} G_{\alpha} (x) = &|2 \sin (\pi x)| \cdot \left| 1+\frac{C-D}{x+1} \right| \cdot \left| 1+\frac{C}{x} \right| \cdot \left| 1+\frac{D}{x-1} \right| \times \\ & \qquad \times \exp \left( C \left( \log \frac{a}{2 \pi} - \frac{\Gamma' (2+x)}{\Gamma (2+x)} \right) +O \left( \frac{1+\log a}{(2-|x|)^2 a^2} \right) \right) \end{split}
\end{equation}
in the range $|x| \le 2-2/a$, where $\Gamma$ is the Gamma function, and
\[ C=\frac{1}{\sqrt{a^2+4}} \qquad \textrm{and} \qquad D=\frac{\sqrt{a^2+4}-a}{2 \sqrt{a^2+4}} . \]
In particular, we roughly have
\begin{equation} \label{Galphalimitfunction2}
G_{\alpha} (x) = |2 \sin (\pi x)| e^{(\log a)/a} +O \left( \frac{1}{a} \right) = |2 \sin (\pi x)| + O \left( \frac{1+\log a}{a} \right) , \quad |x| \le 1.99, 
\end{equation}
but \eqref{Galphalimitfunction} is of course more precise. Observe that the effect of the factors $1+(C-D)/(x+1)$, $1+C/x$, $1+D/(x-1)$ is that they shift the zeroes $-1$, $0$, $1$ of $|2 \sin (\pi x)|$ by roughly $(C-D) \sim 1/a$, $C \sim 1/a$, $D \sim 1/a^2$ to the left, respectively. The admissible range of $x$ in the approximations  \eqref{Galphalimitfunction} and \eqref{Galphalimitfunction2} could be extended by the inclusion of more correction factors; however, in the context of Sudler products only shifts $x$ in the range $x \in (-1,1)$ can occur, so from our perspective there is no reason to aim at a wider range for $x$.

\begin{figure}[h!]%
\centering
\subfigure[][]{%
\label{fig:ex3-a}%
\includegraphics[scale=0.84]{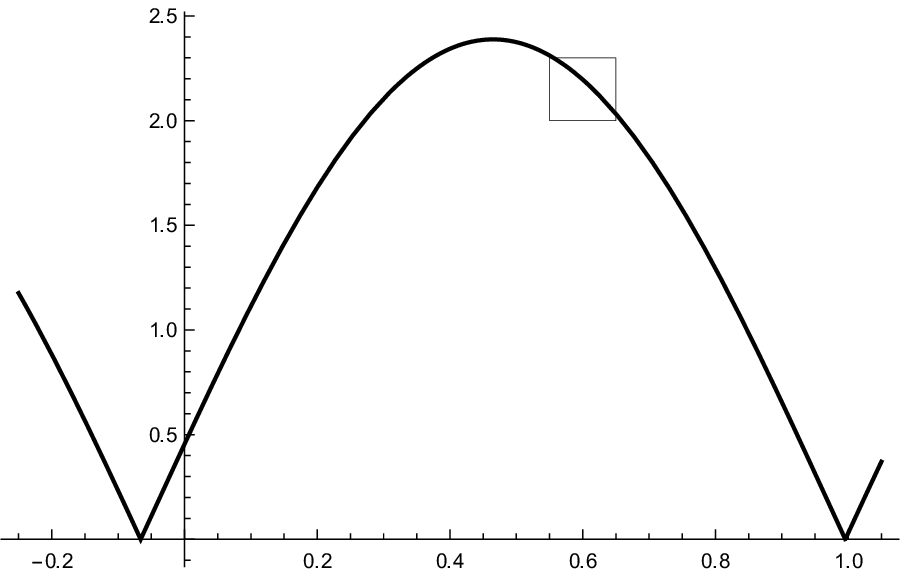}}%
\hspace{8pt}%
\subfigure[][]{%
\label{fig:ex3-b}%
\includegraphics[scale=0.84]{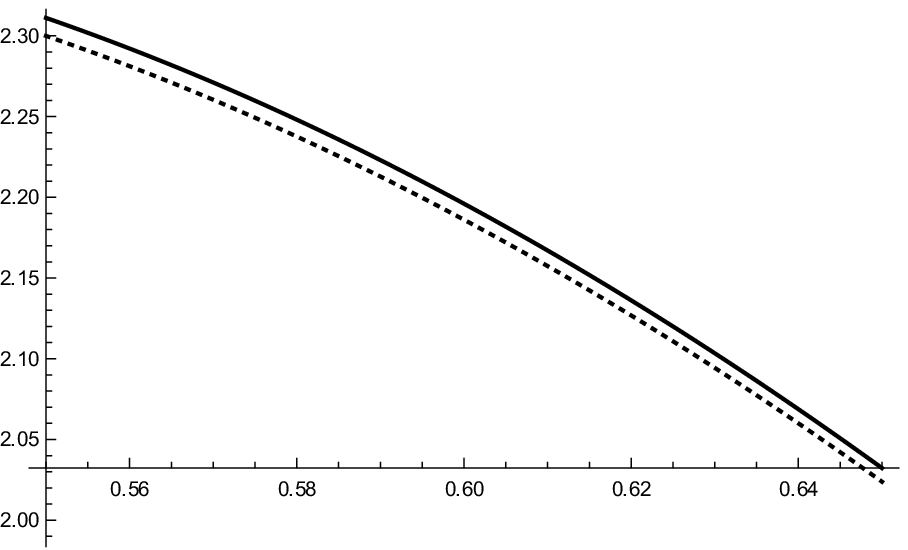}}
\caption{\subref{fig:ex3-a} The function $P_{q_k} (\alpha,(-1)^k x/q_k)$ for $k=4$ and $\alpha = [0;\overline{15}]$, which virtually equals the corresponding limit function $G_\alpha(x)$. \subref{fig:ex3-b} We obtain an excellent approximation from the right-hand side of \eqref{Galphalimitfunction}, leaving out the O-term. The difference is so small that it would be invisible on a full-scale plot as in \subref{fig:ex3-a}, so we have zoomed into the small box indicated in \subref{fig:ex3-a} to show the deviation between the two functions. The actual value of $P_{q_k} (\alpha,x/q_k)$ is plotted as a solid line, the approximation from \eqref{Galphalimitfunction} as a dotted line. Obviously we obtain a much better approximation than the crude $P_{q_k} (\alpha,(-1)^k x/q_k) \approx 2 |\sin(\pi x)|$ of Figure \ref{fig:limf}.}%
\label{fig:ex3}%
\end{figure}

In a recent paper \cite{GNZ} the remarkable convergence property of $P_{q_k} (\alpha, (-1)^k x/q_k)$ was generalized to arbitrary quadratic irrationals $\alpha=[a_0;a_1, \dots, a_{k_0}, \overline{a_{k_0+1}, \dots, a_{k_0+p}}]$; we recall that the overline denotes period. The only difference is that in general we have $p$ different limit functions $G_{\alpha, r} (x)$, $1 \le r \le p$, and $P_{q_k}(\alpha, (-1)^k x/q_k) \to G_{\alpha, r}(x)$ holds pointwise on $\mathbb{R}$ as $k \to \infty$ along the arithmetic progression $k \in p\mathbb{N} + k_0 + r$. Generalizing \eqref{Galphalimitfunction}, the following result states that all these limit functions are close to $|2 \sin (\pi x)|$ whenever the partial quotients of $\alpha$ are all large, and are roughly of similar order of magnitude; the latter property is measured by the parameter $T$.
\begin{thm}\label{quadraticlimitfunction} Let $\alpha=[a_0;a_1, \dots, a_{k_0}, \overline{a_{k_0+1}, \dots, a_{k_0+p}}]$ be a quadratic irrational, and assume that $\max_{1 \le r \le p} (\log a_{k_0+r})/a_{k_0+r+1} \le T$ with some constant $T \ge 1$. For any $1\le r \le p$ and any $|x| \le \max \{ 1, 2-2/a_{k_0+r+1} \}$,
\[ \begin{split} G_{\alpha, r}(x) = &|2 \sin (\pi x)| \cdot \left| 1+\frac{C_r-D_r}{x+1} \right| \cdot \left| 1+\frac{C_r}{x} \right| \cdot \left| 1+\frac{D_r}{x-1} \right| \times \\ &  \times \exp \left( C_r \left( \log \frac{a_{k_0+r}}{2 \pi} - \frac{\Gamma' (2+x)}{\Gamma (2+x)} \right)  \right. \\
& \qquad \qquad \left. + O \left( \frac{T+\log (a_{k_0+r-1}a_{k_0+r})}{(2-|x|)a_{k_0+r} a_{k_0+r+1}} + \frac{T}{(2-|x|)^2 a_{k_0+r+1}^2} \right) \right) , \end{split} \]
where $\Gamma$ is the Gamma function, $C_r=\lim_{m \to \infty} q_{k_0+r+mp} \| q_{k_0+r+mp} \alpha \|$ and $D_r=\lim_{m \to \infty} q_{k_0+r-1+mp} \| q_{k_0+r+mp} \alpha \|$.
\end{thm}

Based on these results for quadratic irrationals with large partial quotients, it is not difficult to come up with the intuition that for a well approximable irrational the corresponding limit function is precisely $|2 \sin (\pi x)|$.
\begin{thm}\label{wellapproximablelimitfunction} Assume that $\sup_{k \ge 1} a_k = \infty$. Then
\[ P_{q_{k_m}} (\alpha, (-1)^{k_m} x/q_{k_m}) \to |2 \sin (\pi x)| \qquad \textrm{as } m \to \infty \]
locally uniformly on $\mathbb{R}$ for any increasing sequence of positive integers $k_m$ such that
\[ \frac{1+\log \max_{1 \le \ell \le k_m} a_{\ell}}{a_{k_m+1}} \to 0 \qquad \textrm{as } m \to \infty . \]
If in addition $\lim_{k \to \infty} (1+\log a_k)/a_{k+1}=0$, then the same holds along the full sequence $k_m=m$.
\end{thm}

\section{Approximation of shifted Sudler products} \label{sec_approx_shift}

The main result of this section is Proposition \ref{pqkapproximation} in Section \ref{subsec_app} below, which is an approximation formula for the inner product over $0 \le b \le b_k-1$ in the decomposition formula \eqref{pnproductform}. As we will see, lower estimates are much more difficult to prove than upper estimates, especially when the points $bq_k \| q_k \alpha \| + \varepsilon_k (N)$ are close to $0$ or $1$, requiring a somewhat tedious case analysis throughout the paper. This is explained by the fact that $\log |2 \sin (\pi x)|$ is bounded above but not below, and has singularities at $x=0$ and $x=1$.

\subsection{Continued fractions} \label{sec_cf}

We start by recalling some basic facts about continued fractions; see \cite{ASA,RSC,SD} for background. The convergents satisfy the recursion $q_{k+1}=a_{k+1} q_k + q_{k-1}$ with initial conditions $q_0=1$, $q_1=a_1$, and $p_{k+1}=a_{k+1}p_k+p_{k-1}$ with initial conditions $p_0=a_0$, $p_1=a_0a_1+1$. If either $k \ge 1$, or $k=0$ and $a_1 >1$, then the following hold:
\begin{enumerate}
\item[(i)] by the best rational approximation property, $\| n \alpha \| \ge \| q_k \alpha \|$ for all $1 \le n < q_{k+1}$;
\item[(ii)] the integer closest to $q_k \alpha$ is $p_k$, and $(-1)^k (q_k \alpha - p_k) = \| q_k \alpha \|$;
\item[(iii)] $1/(q_k \| q_k \alpha \|) = [a_{k+1};a_{k+2},a_{k+3},\cdots ] + [0;a_k,a_{k-1},\dots, a_1]$.
\end{enumerate}
Note that (iii) follows easily from the well-known algebraic identity
\[ [a_0;a_1,\dots, a_k,x]=\frac{p_k x+p_{k-1}}{q_k x+q_{k-1}} \]
with $x=[a_{k+1};a_{k+2},a_{k+3},\dots ]$, and the fact that $q_{k-1}/q_k=[0;a_k,a_{k-1},\dots, a_1]$. In particular, (iii) implies that $1/(a_{k+1}+2) \le q_k \| q_k \alpha \| \le 1/a_{k+1}$.

The recursion $\| q_{k+1} \alpha \| = -a_{k+1} \| q_k \alpha \| + \| q_{k-1} \alpha \|$ and the identity $|\alpha -p_k/q_k| + |\alpha -p_{k+1}/q_{k+1}| =1/(q_k q_{k+1})$, in other words $q_{k+1} \| q_k \alpha \| + q_k \| q_{k+1} \alpha \| =1$, are also classical. Finally, recall the identity
\begin{equation}\label{qk+1pkidentity}
q_{k+1} p_k - q_k p_{k+1} = (-1)^{k+1}, \qquad k \ge 0.
\end{equation}
The Ostrowski expansion of a non-negative integer $N$ is the unique representation $N=\sum_{k=0}^{K-1}b_k q_k$, where $0 \le b_0 <a_1$ and $0 \le b_k \le a_{k+1}$ are integers which satisfy the rule that $b_{k-1}=0$ whenever $b_k=a_{k+1}$.

We first prove a useful estimate for $\varepsilon_k (N)$, as defined in \eqref{varepsilonkdef}. Note that in the product formula \eqref{pnproductform} only those indices $k$ appear for which $b_k \ge 1$; otherwise the inner product is empty, and by convention equals $1$. For all intents and purposes, $\varepsilon_k (N)$ is thus only defined for those $k$ for which $b_k \ge 1$.
\begin{lem}\label{epsilonklemma} Let $N=\sum_{k=0}^{K-1} b_k q_k$ be the Ostrowski expansion of a non-negative integer. For any $k \ge 0$ such that $b_k \ge 1$,
\begin{equation}\label{varepsilonkestimate}
-1 < -q_k \| q_k \alpha \| + q_k \| q_{k+1} \alpha \| \le \varepsilon_k (N) \le q_k \| q_{k+1} \alpha \| < \frac{1}{2} .
\end{equation}
If $b_{k+1} \le (1-\delta ) a_{k+2}$ with some $\delta >0$, then $\varepsilon_k (N) \ge -(1-\delta /3) q_k \| q_k \alpha \|$. If condition \eqref{logak/ak+1} holds, then $\varepsilon_k (N) \ge -(1-\frac{1}{e^T+2})$ for any $k \ge k_0$ such that $b_k \ge 1$.
\end{lem}

\begin{proof} The estimate \eqref{varepsilonkestimate} was already observed in \cite[Lemma 3]{AB}. To see the second claim, assume that $b_{k+1} \le (1-\delta) a_{k+2}$. Then
\[ \begin{split} \varepsilon_k (N) &\ge -q_k \left( b_{k+1} \| q_{k+1} \alpha \| + b_{k+3} \| q_{k+3} \alpha \| + \cdots  \right) \\ &\ge -q_k \left( (1-\delta ) a_{k+2} \| q_{k+1} \alpha \| + a_{k+4} \| q_{k+3} \alpha \| + \cdots \right) \\ &= -q_k \left( (1-\delta ) (\| q_k \alpha \| - \| q_{k+2} \alpha \|) + (\| q_{k+2} \alpha) \| - \| q_{k+4} \alpha \|) + \cdots \right) \\ &= -q_k \left( (1-\delta ) \| q_k \alpha \| + \delta \| q_{k+2} \alpha \| \right) . \end{split} \]
It is not difficult to see that $\| q_{k+2} \alpha \| \le (2/3) \| q_k \alpha \|$. In particular, we have $\varepsilon_k (N) \ge -(1-\delta /3) q_k \| q_k \alpha \|$, as claimed.

To see the last claim, assume that \eqref{logak/ak+1} holds. We then have $q_k \| q_k \alpha \| \le 1-\frac{1}{e^T+2}$ for all $k \ge k_0$. Indeed, this trivially follows from $q_k \| q_k \alpha \| \le 1/a_{k+1}$ if $a_{k+1} \ge 2$. If $a_{k+1}=1$, then property (iii) of continued fractions above gives the more precise bound $1/(q_k \| q_k \alpha \|) \ge 1+1/(a_k+1)$. By condition \eqref{logak/ak+1} here $a_k \le e^T$, and $q_k \| q_k \alpha \| \le 1-\frac{1}{e^T+2}$ follows. Formula \eqref{varepsilonkestimate} thus gives $\varepsilon_k(N) \ge -q_k \| q_k \alpha \| \ge -(1-\frac{1}{e^T+2})$, as claimed.
\end{proof}

\subsection{A cotangent sum}\label{cotangentsection}

The cotangent sum
\begin{equation}\label{cotangentsum}
\sum_{n=1}^{q_k-1} \frac{n}{q_k} \cot \left( \pi \frac{n p_k}{q_k} \right)
\end{equation}
will play an important role in our estimates for the shifted Sudler products. This sum is called ``Vasyunin sum'' after Vasyunin's foundational work in \cite{VA}. It is related to the B\'aez-Duarte--Nyman--Beurling criterion for the Riemann hypothesis; see in particular \cite{MR}. As we already observed in \cite{AB}, a general result of Lubinsky \cite[Theorem 4.1]{LU} implies that for an arbitrary irrational $\alpha$,
\begin{equation}\label{Lubinskybound}
\left| \sum_{n=1}^{q_k-1} \frac{n}{q_k} \cot \left( \pi \frac{n p_k}{q_k} \right) \right| \ll ( 1+\log \max_{1 \le \ell \le k} a_{\ell} ) q_k .
\end{equation}
A reciprocity formula of Bettin and Conrey \cite{BC2} provides a precise evaluation of \eqref{cotangentsum}. In particular, under the assumption \eqref{logak/ak+1} we can isolate a main term; this main term is responsible for the exponential correction factor in Theorem \ref{quadraticlimitfunction}. We now give an approximate evaluation of a shifted version of \eqref{cotangentsum}.
\begin{lem}\label{cotangentevaluation} Assume \eqref{logak/ak+1}. For any $k \ge 4$ and any $x \in (-1,1)$,
\[ \begin{split} \sum_{n=1}^{q_k-1} \frac{n}{q_k} \cot &\left( \pi \frac{n p_k+(-1)^k x}{q_k} \right) \\ &=\frac{(-1)^k q_k}{\pi} \left( \log \frac{a_k}{2 \pi} - \frac{\Gamma' (1+x)}{\Gamma (1+x)} +O \left( \frac{T +\log (a_{k-1}a_k)}{(1-|x|)a_k} \right) \right) + O_{\alpha} (1) , \end{split} \]
where $\Gamma$ is the Gamma function.
\end{lem}

\begin{proof} Let
\[ C_k(x)= \sum_{n=1}^{q_k-1} \frac{n}{q_k} \cot \left( \pi \frac{n p_k+(-1)^k x}{q_k} \right) \]
denote the shifted cotangent sum in the statement of the lemma. We first prove the claim for $x=0$, and then extend it to $x \in (-1,1)$.

It follows from the identity \eqref{qk+1pkidentity} that the multiplicative inverse of $(-1)^{k+1} p_k$ modulo $q_k$ is $q_{k-1}$. We also have the continued fraction expansion $q_{k-1}/q_k = [0;a_k, a_{k-1}, \dots , a_1]$. By a reciprocity formula of Bettin and Conrey \cite{BC2} (see also \cite[Proposition 1]{B}), we have
\begin{equation}\label{BCevaluation}
\begin{split} \sum_{n=1}^{q_k-1} \frac{n}{q_k} \cot \left( \pi \frac{np_k}{q_k} \right) &= (-1)^{k+1} \sum_{n=1}^{q_k-1} \frac{n}{q_k} \cot \left( \pi \frac{n \overline{q_{k-1}}}{q_k} \right)  \\ &= (-1)^{k+1} q_k \sum_{\ell =1}^k \frac{(-1)^{\ell}}{v_{\ell}} \left( \frac{1}{\pi v_{\ell}} + \psi \left( \frac{v_{\ell -1}}{v_{\ell}} \right) \right) , \end{split}
\end{equation}
where $\overline{q_{k-1}}$ denotes the multiplicative inverse of $q_{k-1}$ modulo $q_k$, the fractions $u_{\ell}/v_{\ell} = [0;a_k, a_{k-1}, \dots, a_{k-\ell +1}]$ are the convergents of $q_{k-1}/q_k$ (with the convention $v_0=1$), and $\psi : \mathbb{C} \backslash (-\infty , 0] \to \mathbb{C}$ is an analytic function with asymptotics
\[ \begin{split} \psi (x) = \frac{\log \frac{1}{2 \pi x} + \gamma}{\pi x} + O (\log (1/x)) \end{split} \]
as $x \to 0$ along the positive reals, with $\gamma$ denoting the Euler--Mascheroni constant. The $\ell =1$ term is
\[ \begin{split} \frac{-1}{v_1} \left( \frac{1}{ \pi v_1} + \psi \left( \frac{v_0}{v_1} \right) \right) &= \frac{-1}{a_k} \left( \frac{\log \frac{a_k}{2 \pi} + \gamma}{\pi /a_k} + O(1+ \log a_k) \right) \\ &= \frac{-1}{\pi} \left( \log \frac{a_k}{2 \pi} + \gamma + O \left( \frac{1+\log a_k}{a_k} \right) \right) . \end{split} \]
The terms $2 \le \ell \le k-k_0$ are negligible due to the assumption $(\log a_k)/a_{k+1} \le T$:
\[ \begin{split} \sum_{\ell =2}^{k-k_0} \frac{1}{v_{\ell}} \left| \frac{1}{\pi v_{\ell}} + \psi \left( \frac{v_{\ell -1}}{v_{\ell}} \right) \right| &\ll \sum_{\ell =2}^{k-k_0} \frac{1}{v_{\ell}} \cdot \frac{1+\log (v_{\ell}/v_{\ell -1})}{v_{\ell -1}/v_{\ell}} \\ &\ll \sum_{\ell =2}^{k-k_0} \frac{1+\log a_{k-\ell +1}}{v_{\ell -1}} \\ &\ll \frac{1+\log a_{k-1}}{a_k} + \sum_{\ell =3}^{k-k_0} \frac{T}{v_{\ell -2}} \\ &\ll \frac{1+\log a_{k-1}}{a_k} + \frac{T}{a_k} . \end{split} \]
Finally, the terms $k-k_0+1 \le \ell \le k$ satisfy
\[ \sum_{\ell =k-k_0+1}^k \frac{1}{v_{\ell}} \left| \frac{1}{\pi v_{\ell}} + \psi \left( \frac{v_{\ell -1}}{v_{\ell}} \right) \right| \ll \sum_{\ell =k-k_0+1}^k \frac{1+\log a_{k-\ell +1}}{v_{\ell -1}} \ll_{\alpha} \frac{1}{v_k} = \frac{1}{q_k} . \]
Using the previous three formulas in \eqref{BCevaluation}, we get
\begin{equation}\label{cotangentx=0}
C_k(0) =\frac{(-1)^k q_k}{\pi} \left( \log \frac{a_k}{2 \pi} + \gamma +O \left( \frac{T+\log (a_{k-1}a_k)}{a_k} \right) \right) + O_{\alpha} (1) .
\end{equation}
This proves the claim when $x=0$; note that $-\Gamma' (1)/\Gamma (1) = \gamma$.

Next, let $x \in (-1,1)$, and consider the derivative
\begin{align}
C_k'(x) & = \sum_{n=1}^{q_k-1} \frac{(-1)^{k+1}\pi n}{q_k^2\sin^2 (\pi ((np_k+(-1)^k x)/q_k))} \nonumber \\
& = \sum_{n=1}^{q_k-1} \frac{(-1)^{k+1}n}{\pi q_k^2 \| (n(-1)^k p_k+x)/q_k \|^2} + O(1).\label{cotangentderivative}
\end{align}
In the second step we used the general estimate $\pi/\sin^2 (\pi y) = 1/(\pi \| y \|^2) + O(1)$. We now isolate a small number of integers $n$ which give the main contribution in \eqref{cotangentderivative}. Recall once again that $q_{k-1}p_k \equiv (-1)^{k+1} \pmod{q_k}$. Let $0<|a| \le a_k$ be an integer. Then the solution of the congruence $n(-1)^k p_k \equiv a \pmod{q_k}$ is $n \equiv -a q_{k-1} \pmod{q_k}$; the unique representative of this residue class in $1 \le n \le q_k-1$ is $n=q_k-aq_{k-1}$ if $1 \le a \le a_k$, and $n=-aq_{k-1}$ if $-a_k \le a \le -1$. The contribution of these $2a_k$ integers $n$ in \eqref{cotangentderivative} is
\[ \begin{split} \frac{(-1)^{k+1}}{\pi} \bigg( \sum_{a=1}^{a_k} &\frac{q_k-aq_{k-1}}{q_k^2 \| a/q_k+x/q_k \|^2} + \sum_{a=-a_k}^{-1} \frac{-aq_{k-1}}{q_k^2 \| a/q_k + x/q_k \|^2}  \bigg) \\ &= \frac{(-1)^{k+1}}{\pi} \left( \sum_{a=1}^{a_k} \frac{q_k}{(a+x)^2} + \sum_{a=1}^{a_k} \left( \frac{a q_{k-1}}{(a-x)^2} - \frac{a q_{k-1}}{(a+x)^2} \right) \right) \\ &= \frac{(-1)^{k+1} q_k}{\pi} \sum_{a=1}^{\infty} \frac{1}{(a+x)^2} + O \left( \frac{q_k(1+\log a_k)}{(1-|x|)^2 a_k} \right) . \end{split} \]
Note that we used the assumption $k \ge 4$ to ensure that $(a_k+1)/q_k \le 1/2$. Since the contribution of all other integers $n$ in \eqref{cotangentderivative} is
\[ \ll \sum_{a_k<|a| \le q_k/2} \frac{1}{q_k(1-|x|)^2\| a/q_k \|^2} \ll \frac{q_k}{(1-|x|)^2 a_k} , \]
we get
\[ C_k'(x) = \frac{(-1)^{k+1} q_k}{\pi} \sum_{a=1}^{\infty} \frac{1}{(a+x)^2} + O \left( \frac{q_k (1+\log a_k)}{(1-|x|)^2 a_k} \right) . \]
By integrating and identifying the resulting infinite series as a special function we get
\[ \begin{split} C_k(x)-C_k(0) &= \frac{(-1)^{k+1} q_k}{\pi} \sum_{a=1}^{\infty} \left( \frac{1}{a} - \frac{1}{a+x} \right) + O \left( \frac{q_k (1+\log a_k)}{(1-|x|)a_k} \right) \\ &= \frac{(-1)^{k+1} q_k}{\pi} \left( \gamma + \frac{\Gamma' (1+x)}{\Gamma (1+x)} + O \left( \frac{1+\log a_k}{(1-|x|)a_k} \right) \right) , \end{split} \]
and the claim for general $x \in (-1,1)$ follows from the special case \eqref{cotangentx=0}.
\end{proof}

\subsection{A modified cotangent sum}\label{modifiedcotangentsection}

We will actually need a slightly modified version of the cotangent sum in Lemma \ref{cotangentevaluation}, defined as
\begin{equation}\label{Vkxdef}
V_k (x):= \sum_{n=1}^{q_k-1} \sin (\pi n \| q_k \alpha \| /q_k) \cot \left( \pi \frac{n(-1)^k p_k+x}{q_k} \right) .
\end{equation}

\begin{lem}\label{Vkxlemma}
\leavevmode
\begin{enumerate}
\item[(i)] For any $k \ge 1$, the derivative of $V_k(x)$ on the interval $(-1,1)$ satisfies
\[ V_k'(x) < 0 \quad \textrm{and} \quad |V_k'(x)| \ll \frac{1}{(1-|x|)^2 a_{k+1}} . \] 
\item[(ii)] For any $k \ge 1$,
\[ |V_k(0)| \ll \frac{1+\log \max_{1 \le \ell \le k} a_{\ell}}{a_{k+1}} . \]
\item[(iii)] Assume \eqref{logak/ak+1}. For any $k \ge 4$ and any $x \in (-1,1)$,
\[ \frac{V_k(x)}{q_k \| q_k \alpha \|} = \log \frac{a_k}{2 \pi} - \frac{\Gamma' (1+x)}{\Gamma (1+x)}+ O \left( \frac{T+ \log (a_{k-1} a_k)}{(1-|x|)a_k} \right) + O_{\alpha} \left( \frac{1}{q_k} \right) . \]
In particular,
\[ V_k(x) = \frac{\log a_k}{a_{k+1}} + O \left( \frac{T}{(1-|x|)a_{k+1}} \right) + O_{\alpha} \left( \frac{1}{q_{k+1}} \right) . \]
\end{enumerate}
\end{lem}

\begin{proof} Let $x \in (-1,1)$. Clearly,
\[ V_k'(x) = \sum_{n=1}^{q_k-1} \sin (\pi n \| q_k \alpha \| /q_k) \frac{-\pi}{q_k \sin^2 \left( \pi (n (-1)^k p_k/q_k + x/q_k) \right)} <0 . \]
By the general inequality $|\sin (\pi y)| \ge 2 \| y \|$ and
\[ \| n (-1)^k p_k/q_k + x/q_k  \| \ge (1-|x|) \| np_k/q_k \|, \]
we also have
\[ |V_k'(x)| \ll \sum_{n=1}^{q_k-1} \frac{ \| q_k \alpha \|}{q_k (1-|x|)^2 \| n p_k /q_k \|^2} \ll \frac{q_k  \| q_k \alpha \|}{(1-|x|)^2} \ll \frac{1}{(1-|x|)^2 a_{k+1}} . \]
In the second step we used the fact that as $n$ runs in the interval $1 \le n \le q_k-1$, the integers $np_k$ attain each nonzero residue class modulo $q_k$ exactly once. This finishes the proof of (i).

Next, note that the general estimate $\sin y = y + O(|y|^3)$ implies that the error of replacing $\sin (\pi n \| q_k \alpha \| / q_k)$ by $\pi n \| q_k \alpha \| /q_k$ in the definition of $V_k(x)$ is
\[ \ll \sum_{n=1}^{q_k-1} \frac{n^3 \| q_k \alpha \|^3}{q_k^3} \left| \cot \left( \pi \frac{n (-1)^k p_k+x}{q_k} \right) \right| \ll \sum_{n=1}^{q_k -1} \frac{\| q_k \alpha \|^3}{(1-|x|) \| np_k /q_k \|} \ll \frac{\| q_k \alpha \|^3 q_k \log q_k}{1-|x|} , \]
hence
\[ V_k(x) = \pi \| q_k \alpha \| \sum_{n=1}^{q_k-1} \frac{n}{q_k} \cot \left( \pi \frac{n(-1)^k p_k+x}{q_k} \right) + O \left( \frac{\| q_k \alpha \|^3 q_k \log q_k}{1-|x|} \right) . \]
Claims (ii) and (iii) thus follow from \eqref{Lubinskybound} and Lemma \ref{cotangentevaluation}, respectively.
\end{proof}

\subsection{The reflection and transfer principles}

In our previous paper \cite{AB} we showed the useful identity
\[ P_N (p/q) \cdot P_{q-N-1} (p/q) = q \]
for any reduced fraction $p/q$ and any integer $0 \le N < q$. We also proved that
\[ \left| \log P_N (p_k/q_k) - \log P_N (\alpha) \right| \ll \frac{1+\log \max_{1 \le \ell \le k} a_{\ell}}{a_{k+1}} \]
for an arbitrary irrational $\alpha$ and all $0 \le N < q_k$. We called these results the reflection and transfer principles, respectively; the latter terminology comes from the fact that it helps transfer results between rational and irrational settings. In this section we establish similar principles for shifted Sudler products.
\begin{prop}\label{reflectionprop} Let $p/q$ be a reduced fraction. For any $0 \le N <q$ and any $x \in \mathbb{R}$,
\begin{equation}\label{reflectionprinciple}
P_N (p/q,x) \cdot P_{q-N-1} (p/q,-x) = \left\{ \begin{array}{ll} \frac{|\sin (\pi qx)|}{|\sin (\pi x)|} & \textrm{if } x \not\in \mathbb{Z}, \\ q & \textrm{if } x \in \mathbb{Z} . \end{array} \right.
\end{equation}
In particular, for any $x \in \mathbb{R}$,
\begin{equation}\label{lastterm}
P_{q-1} (p/q,x) = \left\{ \begin{array}{ll} \frac{|\sin (\pi qx)|}{|\sin (\pi x)|} & \textrm{if } x \not\in \mathbb{Z}, \\ q & \textrm{if } x \in \mathbb{Z} . \end{array} \right.
\end{equation}
\end{prop}

\begin{proof} For a given $x \in \mathbb{R}$ consider the factorization
\[ t^q - e^{2 \pi i q x} = (t-e^{2 \pi i x}) \prod_{j=1}^{q-1} (t-e^{2 \pi i (j/q+x)}) . \]
Dividing both sides by $(t-e^{2 \pi i x})$ and letting $t \to 1$, we get
\[  \prod_{j=1}^{q-1} (1-e^{2 \pi i (j/q+x)}) = \left\{ \begin{array}{ll} \frac{1-e^{2 \pi i q x}}{1-e^{2 \pi i x}} & \textrm{if } x \not\in \mathbb{Z}, \\ q & \textrm{if } x \in \mathbb{Z}.  \end{array} \right.  \]
Therefore
\[ P_{q-1} (p/q,x) = \prod_{n=1}^{q-1} |1-e^{2 \pi i (np/q+x)}| = \prod_{j=1}^{q-1} |1-e^{2 \pi i (j/q+x)}| = \left\{ \begin{array}{ll} \frac{|\sin (\pi qx)|}{|\sin (\pi x)|} & \textrm{if } x \not\in \mathbb{Z}, \\ q & \textrm{if } x \in \mathbb{Z} , \end{array} \right. \]
as claimed in \eqref{lastterm}. Next, let $0 \le N < q$. By the definition of shifted Sudler products and the previous formula,
\[ P_N (p/q,x) \cdot \prod_{n=N+1}^{q-1} |2 \sin (\pi (np/q +x))| = P_{q-1} (p/q, x) = \left\{ \begin{array}{ll} \frac{|\sin (\pi qx)|}{|\sin (\pi x)|} & \textrm{if } x \not\in \mathbb{Z}, \\ q & \textrm{if } x \in \mathbb{Z} . \end{array} \right. \]
A simple reindexing shows that here
\[ \prod_{n=N+1}^{q-1} |2 \sin (\pi (np/q +x))| = \prod_{j=1}^{q-N-1} |2 \sin (\pi ((q-j)p/q+x))| = P_{q-N-1} (p/q, -x) , \]
which proves \eqref{reflectionprinciple}.
\end{proof}

\begin{cor}\label{Bqk-1corollary} Let $k \ge 1$ and $0 \le M <q_k$ be integers, and define
\begin{equation}\label{BNx}
B_{k,M}(x):= \log \frac{P_M (\alpha, (-1)^k x/q_k)}{P_M (p_k/q_k, (-1)^k x/q_k)} - \sum_{n=1}^M \sin (\pi n \| q_k \alpha \| /q_k ) \cot \left( \pi \frac{n (-1)^k p_k+x}{q_k} \right) .
\end{equation}
Then
\[ \log P_{q_k}(\alpha, (-1)^k x/q_k) = \log \left( |2 \sin (\pi (\| q_k \alpha \| + x/q_k))| \frac{|\sin (\pi x)|}{|\sin (\pi x/q_k)|} \right) + V_k (x) + B_{k,q_k-1}(x) , \]
with the convention $|\sin (\pi x)| / |\sin (\pi x /q_k)| =q_k$ when $x/q_k \in \mathbb{Z}$.
\end{cor}

\begin{proof} By the definitions \eqref{BNx} of $B_{k,M}(x)$ and \eqref{Vkxdef} of $V_k(x)$,
\begin{equation}\label{transferqk-1}
\log P_{q_k-1} (\alpha, (-1)^k x /q_k) = \log P_{q_k-1} (p_k/q_k, (-1)^k x /q_k) +V_k(x) + B_{k,q_k-1}(x) .
\end{equation}
Using the identity \eqref{lastterm}, here
\[ \log P_{q_k-1} (p_k/q_k, (-1)^k x /q_k) = \log \frac{|\sin (\pi x)|}{|\sin (\pi x/q_k)|} . \]
Adding $\log |2 \sin (\pi (q_k \alpha + (-1)^k x/q_k))| = \log |2 \sin (\pi ( \| q_k \alpha \| + x/q_k))|$ to both sides of \eqref{transferqk-1}, the claim follows.
\end{proof}

\noindent In the claim of Corollary \ref{Bqk-1corollary} we consider $V_k(x)$ to be a first order correction term, and $B_{k,q_k-1}(x)$ to be an error term. The following proposition gives estimates for $B_{k,M}(x)$; we call it the transfer principle for shifted Sudler products. In fact, in the present paper we will only use it with $M=q_k-1$.
\begin{prop}\label{transferprop}
\leavevmode
\begin{enumerate}
\item[(i)] Let $k \ge 1$ and $0 \le M < q_k$ be integers, and assume that $q_k \| q_k \alpha \| \le 1-c_k$ and $-1<x \le 1-\frac{q_k \| q_k \alpha \|}{1-c_k}$ for some $c_k$ such that $10/q_k^2 \le c_k<1$. Then
\[ -C \frac{\log (4/c_k)}{(1-|x|)^2 a_{k+1}^2} \le B_{k,M}(x) \le C \frac{1}{a_{k+1}^2 q_k} \]
with a universal constant $C>0$.
\item[(ii)] Let $N=\sum_{k=0}^{K-1} b_k q_k$ be the Ostrowski expansion of a non-negative integer. For any $1 \le k \le K-1$, any $0 \le M <q_k$ and any $0 \le b \le b_k-1$, we have
\[ B_{k,M}(b q_k \| q_k \alpha \| + \varepsilon_k (N)) \le C \frac{1}{a_{k+1}^2q_k} \]
with a universal constant $C>0$.
\end{enumerate}
\end{prop}

\begin{proof}[Proof of Proposition \ref{transferprop} (i)] Using trigonometric identities we can write
\begin{equation}\label{prod1+xnshifted}
\frac{P_M (\alpha, (-1)^k x/q_k )}{P_M (p_k/q_k, (-1)^k x/q_k )} = \left| \prod_{n=1}^M \frac{\sin (\pi (n \alpha +(-1)^k x/q_k))}{\sin (\pi (n p_k/q_k +(-1)^k x/q_k))} \right| = \left| \prod_{n=1}^M (1+x_n+y_n) \right|,
\end{equation}
where
\[ x_n := \cos (\pi n (\alpha -p_k/q_k )) -1= \cos (\pi n \| q_k \alpha \| /q_k) -1 \]
and
\[ \begin{split} y_n &:= \sin (\pi n (\alpha -p_k/q_k)) \cot (\pi (n p_k/q_k +(-1)^k x/q_k)) \\ &= \sin (\pi n \| q_k \alpha \| /q_k ) \cot (\pi (n (-1)^k p_k/q_k + x/q_k)) . \end{split} \]
Assume first, that $0 \le x \le 1-\frac{q_k\| q_k\alpha \|}{1-c_k}$. From the Taylor expansions of sine and cosine, and the estimate
\begin{equation}\label{npkqk+xqk}
\| n(-1)^k p_k/q_k+x/q_k \| \ge (1-x) \| np_k /q_k \| \ge (1-x)/q_k
\end{equation}
we get that for any $0<n<q_k$,
\begin{equation}\label{cosineboundshifted}
|x_n| \le \frac{\pi^2 n^2 \| q_k \alpha \|^2}{2q_k^2} \le \frac{c_k}{2}
\end{equation}
and
\[ \begin{split} \left| y_n \right| \le \frac{\sin (\pi n \| q_k \alpha \| /q_k)}{|\sin (\pi (n(-1)^k p_k/q_k + x/q_k))|} &\le \frac{\pi \| q_k \alpha \|}{\pi (1-x) /q_k - \pi^3 (1-x)^3 /(6q_k^3)} \\ &\le \frac{q_k \| q_k \alpha \|}{1-x} \cdot \frac{1}{1-\pi^2 /(6q_k^2)} \\ & \le (1-c_k) \frac{1}{1-\pi^2 /(6q_k^2)} \\ &\le 1-\frac{3c_k}{4} . \end{split} \]
The point is that each factor in \eqref{prod1+xnshifted} is bounded away from zero, as $1+x_n+y_n \ge c_k/4$; in particular, the absolute values in \eqref{prod1+xnshifted} can be removed. Since $y_n$ is a decreasing function of $x \in (-1,1)$, the same holds if $-1<x<0$.

Observe that for any $t \ge -1+c_k/4$,
\[ e^{t-2t^2 \log (4/c_k)} \le 1+t \le e^t . \]
Indeed, one readily verifies that the function $e^{-t+2t^2 \log (4/c_k)} (1+t)$ attains its minimum on the interval $[-1+c_k/4, \infty)$ at $t=0$. Applying this estimate with $t=x_n+y_n$ in each factor of \eqref{prod1+xnshifted}, we obtain
\begin{equation}\label{pnalpha/pnpkqk}
\begin{split} \exp \left( \sum_{n=1}^M (x_n+y_n) -2 \sum_{n=1}^M (x_n+y_n)^2 \log (4/c_k) \right) &\le \frac{P_M (\alpha , (-1)^k x/q_k )}{P_M (p_k/q_k, (-1)^k x/q_k )} \\ &\le \exp \left( \sum_{n=1}^M (x_n+y_n) \right) . \end{split}
\end{equation}
By \eqref{cosineboundshifted}, we have
\[ \sum_{n=1}^M |x_n| \ll \sum_{n=1}^M \frac{1}{a_{k+1}^2 q_k^2} \ll \frac{1}{a_{k+1}^2 q_k} \]
and
\[ \sum_{n=1}^M x_n^2 \log (4/c_k) \ll \sum_{n=1}^M \frac{\log (4/c_k)}{a_{k+1}^4 q_k^4} \ll \frac{\log (4/c_k)}{a_{k+1}^4 q_k^3} . \]
From \eqref{npkqk+xqk} we get
\[ |y_n| \le \frac{\sin (\pi n \| q_k \alpha \| / q_k)}{|\sin (\pi (n(-1)^k p_k/q_k +x/q_k))|} \ll \frac{\|q_k \alpha \|}{(1-|x|) \| n p_k /q_k \|} , \]
and hence
\[ \sum_{n=1}^M y_n^2 \log (4/c_k) \ll \sum_{n=1}^M \frac{\| q_k \alpha \|^2 \log (4/c_k)}{(1-|x|)^2 \| np_k/q_k \|^2} \ll \frac{\log (4/c_k)}{(1-|x|)^2 a_{k+1}^2} . \]
The estimate \eqref{pnalpha/pnpkqk} thus simplifies as
\[ -C \frac{\log (4/c_k)}{(1-|x|)^2 a_{k+1}^2} \le \log \frac{P_M (\alpha, (-1)^k x/q_k)}{P_M (p_k/q_k, (-1)^k x/q_k)} - \sum_{n=1}^M y_n \le C \frac{1}{a_{k+1}^2 q_k} \]
with some universal constant $C>0$, which proves the claim.
\end{proof}

\begin{proof}[Proof of Proposition \ref{transferprop} (ii)] We argue as in the previous proof. First, we claim that in \eqref{prod1+xnshifted} the absolute values can be removed at the point $x=b q_k \| q_k \alpha \| + \varepsilon_k (N)$. To see this, note that
\[ n \alpha + (-1)^k x/q_k = (-1)^k ((n+bq_k) \| q_k \alpha \| + \varepsilon_k (N))/q_k +np_k/q_k . \]
By Lemma \ref{epsilonklemma}, here
\[ \begin{split} (n+bq_k) \| q_k \alpha \| + \varepsilon_k (N) &\le (b+1)q_k \| q_k \alpha \| + q_k \| q_{k+1} \alpha \| \\ &< q_{k+1} \| q_k \alpha \| + q_k \| q_{k+1} \alpha \| \\ &=1 , \end{split} \]
and also
\[ (n+bq_k) \| q_k \alpha \| + \varepsilon_k (N) \ge \varepsilon_k (N) > -1 . \]
Consequently, $|n \alpha + (-1)^k x/q_k -np_k/q_k|<1/q_k$. We clearly also have $|x|<1$, therefore the points $n \alpha + (-1)^k x/q_k$ and $np_k/q_k + (-1)^k x/q_k$ both lie in the open interval centered at $np_k/q_k \not\in \mathbb{Z}$ of radius $1/q_k$. Since the function $\sin (\pi y)$ does not have a zero in this interval, we have
\[ \frac{\sin (\pi (n \alpha + (-1)^k x/q_k))}{\sin (\pi (np_k/q_k + (-1)^k x/q_k))} >0 . \]
Hence \eqref{prod1+xnshifted} indeed holds without the absolute values; that is,
\[ \frac{P_M (\alpha, (-1)^k x/q_k )}{P_M (p_k/q_k, (-1)^k x/q_k )} = \prod_{n=1}^M \frac{\sin (\pi (n \alpha +(-1)^k x/q_k))}{\sin (\pi (n p_k/q_k +(-1)^k x/q_k))} = \prod_{n=1}^M (1+x_n+y_n) \]
with $x_n$, $y_n$ as in the previous proof. The upper bound
\[ \frac{P_M (\alpha, (-1)^k x/q_k )}{P_M (p_k/q_k, (-1)^k x/q_k )} \exp \left( -\sum_{n=1}^M y_n \right) \le \exp \left( \sum_{n=1}^M x_n \right) \le \exp \left( C \frac{1}{a_{k+1}^2q_k} \right) \]
immediately follows, as claimed.
\end{proof}

\subsection{Key estimate for shifted Sudler products}\label{subsec_app}

We emphasize that in the following proposition we do not assume condition \eqref{logak/ak+1}, hence it could serve as a starting point for various generalizations of the results in this paper. In the proofs of our theorems, condition \eqref{logak/ak+1} will ensure that in the claim of the proposition the contribution of the cotangent sum (the sum expressed in terms of $V_k(x)$) is negligible compared to the sum which is expressed in terms of $\log |2 \sin (\pi x )|$.

\begin{prop}\label{pqkapproximation} Let $N=\sum_{k=0}^{K-1} b_k q_k$ be the Ostrowski expansion of a non-negative integer. For any $k \ge 1$ such that $b_k \ge 1$,
\[ \begin{split} \sum_{b=0}^{b_k-1} \log P_{q_k} (\alpha, (-1)^k (b q_k \| q_k \alpha \| + \varepsilon_k (N))/q_k) = &\sum_{b=1}^{b_k-1} \log |2 \sin (\pi (b q_k \| q_k \alpha \| + \varepsilon_k (N) ))| \\ &+\sum_{b=0}^{b_k-1} V_k (b q_k \| q_k \alpha \| + \varepsilon_k (N) ) \\ &+ \log (2 \pi (b_k q_k \| q_k \alpha \| + \varepsilon_k (N) )) \\&+ E_k(N) ,  \end{split} \]
where $E_k(N) \le C /(a_{k+1}q_k)$ with a universal constant $C>0$. If in addition $k \ge 20 \log \frac{20}{\delta}$, $b_k \le (1-\delta) a_{k+1}$ and $q_k \| q_k \alpha \| \le 1-\delta$ with some $\delta >0$, then we also have $E_k(N) \ge -C\frac{\log (2/\delta) }{\delta^2}(1/a_{k+1} + 1/q_k^2)$ with a universal constant $C>0$.
\end{prop}

\begin{proof} For the sake of readability, put $f(x)=|2 \sin (\pi x)|$ and $\varepsilon_k=\varepsilon_k (N)$. Applying Corollary \ref{Bqk-1corollary} at $x=b q_k \| q_k \alpha \| + \varepsilon_k$ and summing over $0 \le b \le b_k-1$, we get
\[ \begin{split} \sum_{b=0}^{b_k-1} \log P_{q_k} (\alpha, (-1)^k &(b q_k \| q_k \alpha \| + \varepsilon_k )/q_k) \\ =&\sum_{b=0}^{b_k-1} \log \left( f((b+1)\| q_k \alpha \| + \varepsilon_k /q_k ) \frac{f(b q_k \|q_k \alpha \| + \varepsilon_k )}{f(b \|q_k \alpha \| + \varepsilon_k /q_k)} \right) \\ &+ \sum_{b=0}^{b_k-1} V_k (b q_k \| q_k \alpha \| + \varepsilon_k ) \\ &+ \sum_{b=0}^{b_k-1} B_{k,q_k-1} (b q_k \| q_k \alpha \| + \varepsilon_k ) . \end{split} \]
Observe that the first sum on the right hand side has a telescoping part. By peeling off the $b=0$ term we thus obtain
\[ \begin{split} \sum_{b=0}^{b_k-1} \log P_{q_k} (\alpha, (-1)^k (b q_k \| q_k \alpha \| + \varepsilon_k )/q_k) =&\sum_{b=1}^{b_k-1} \log f( b q_k \| q_k \alpha \| + \varepsilon_k ) \\ &+ \sum_{b=0}^{b_k-1} V_k (b q_k \| q_k \alpha \| + \varepsilon_k ) \\ &+ \log \left( f(b_k \| q_k \alpha \| +\varepsilon_k /q_k ) \frac{f(\varepsilon_k )}{f(\varepsilon_k /q_k)} \right) \\ &+ \sum_{b=0}^{b_k-1} B_{k,q_k-1} (b q_k \| q_k \alpha \| + \varepsilon_k ) , \end{split} \]
with the convention that $f(\varepsilon_k ) / f(\varepsilon_k /q_k) =q_k$ if $\varepsilon_k =0$. It remains to estimate the error term
\[ E_k(N) := \log \left( \frac{f(b_k \| q_k \alpha \| +\varepsilon_k /q_k )}{2 \pi (b_k q_k \| q_k \alpha \| + \varepsilon_k )} \cdot \frac{f(\varepsilon_k )}{f(\varepsilon_k /q_k)} \right) + \sum_{b=0}^{b_k-1} B_{k,q_k-1} (b q_k \| q_k \alpha \| + \varepsilon_k ) . \]

First, we prove the upper bound for $E_k(N)$. Using $B_{k,q_k-1}(bq_k \| q_k \alpha \| +\varepsilon_k ) \le C / (a_{k+1}^2 q_k)$ from Proposition \ref{transferprop} (ii) and elementary estimates for the sine function,
\[ E_k(N) \le \log \left( \frac{f(b_k \| q_k \alpha \| +\varepsilon_k /q_k )}{2 \pi (b_k q_k \| q_k \alpha \| + \varepsilon_k )} \cdot q_k \right) + \sum_{b=0}^{b_k-1} \frac{C}{a_{k+1}^2q_k} \le \frac{C}{a_{k+1}q_k},  \]
as claimed.

Next, assume in addition, that $k \ge 20 \log \frac{20}{\delta}$, $b_k \le (1-\delta ) a_{k+1}$ and $q_k \| q_k \alpha \| \le 1-\delta$. By Lemma \ref{epsilonklemma} for any $0 \le b \le b_k-1$, the point $x=bq_k \| q_k \alpha \| + \varepsilon_k$ satisfies
\[ \begin{split} x &\le ((1-\delta ) a_{k+1}-1) q_k \| q_k \alpha \| + q_k \| q_{k+1} \alpha \| \\ &=((1-\delta ) a_{k+1}-1) q_k \| q_k \alpha \| +1-q_{k+1} \| q_k \alpha \| \\&\le 1-(1+\delta ) q_k \| q_k \alpha \| , \end{split} \]
and also
\[ x \ge \varepsilon_k \ge -q_k \| q_k \alpha \| \ge -(1-\delta ) . \]
Hence we can apply Proposition \ref{transferprop} (i) with $c_k=\delta /(1+\delta )$. Note that the assumption $k \ge 20 \log \frac{20}{\delta}$ ensures that $c_k \ge 10/q_k^2$. Since we also have $|x| \le 1-\delta$, we obtain
\[ \sum_{b=0}^{b_k-1} B_{k,q_k-1} (b q_k \| q_k \alpha \| + \varepsilon_k ) \ge \sum_{b=0}^{b_k-1} \left(- \frac{C \log (4/c_k)}{\delta^2 a_{k+1}^2} \right) \ge - \frac{C \log (2/\delta )}{\delta^2 a_{k+1}} . \]
Finally, using the general estimate $\sin y = y(1+O(y^2))$ we get
\[ f(b_k \| q_k \alpha \| + \varepsilon_k /q_k) = 2 \pi (b_k \|q_k \alpha \| + \varepsilon_k /q_k)  (1+O( 1/q_k^2 )) \]
and
\[ \frac{f(\varepsilon_k )}{f(\varepsilon_k /q_k)} = q_k (1+O(\varepsilon_k^2)) = q_k ( 1+O( 1/a_{k+1}^2 )) . \]
Therefore
\[ \begin{split} \log \left( \frac{f(b_k \| q_k \alpha \| +\varepsilon_k /q_k )}{2 \pi (b_k q_k \| q_k \alpha \| + \varepsilon_k )} \cdot \frac{f(\varepsilon_k )}{f(\varepsilon_k /q_k)} \right) &= \log (1+O(1/a_{k+1}^2 + 1/q_k^2 )) \\ &= O(1/a_{k+1}^2 + 1/q_k^2 ) .  \end{split} \]
Altogether we get $E_k(N) \ge -C \frac{\log (2/\delta)}{\delta^2} (1/a_{k+1} + 1/q_k^2)$, as claimed.
\end{proof}

\section{Proof of Theorem \ref{maxpntheorem}} \label{sec_proofs_1}

Throughout this section we assume that \eqref{logak/ak+1} holds with some $k_0, T \ge 1$. Let $\delta_T >0$ be a small enough constant depending only on $T$; for the convenience of the reader we mention that $\delta_T = \min\{ 1/(4 \pi e^{2T}), 1/100 \}$ is a suitable choice. We may assume that $k_0 \ge 20 \log \frac{20}{\delta_T}$. Let us now introduce a sequence which will play a key role in the proof of Theorem \ref{maxpntheorem}.
\begin{definition}\label{UNdef} Given a non-negative integer with Ostrowski expansion $N=\sum_{k=0}^{K-1} b_k q_k$, for any $k_0 \le k \le K-1$ let $u_k(N)=1$ if $b_k=0$, and let
\[ \begin{split} u_k(N) = &\bigg( \prod_{b=1}^{b_k-1} |2 \sin (\pi (b q_k \| q_k \alpha \| + \varepsilon_k (N)))| \bigg) \exp \left( \sum_{b=0}^{b_k-1} V_k (b q_k \| q_k \alpha \| + \varepsilon_k (N) ) \right) \times \\ & \times 2 \pi (b_k q_k \| q_k \alpha \| + \varepsilon_k (N)) \end{split} \]
if $b_k \ge 1$. Finally, let $U_N = \prod_{k=k_0}^{K-1} u_k(N)$.
\end{definition}
\noindent Summarizing the results of the previous section, we can rephrase Proposition \ref{pqkapproximation} in terms of $U_N$.
\begin{prop}\label{PNUNprop} For any non-negative integer with Ostrowski expansion $N=\sum_{k=0}^{K-1} b_k q_k$, we have
\[ \log P_N (\alpha ) = \log U_N - \sum_{k=k_0}^{K-1} F_k(N) + O_T \left( \sum_{k=1}^K \frac{1}{a_k} \right) + O_{\alpha} (1) \]
with some $F_k(N)$ satisfying $F_k(N) \ge 0$ for all $k_0 \le k \le K-1$, and $F_k(N)=0$ for all $k_0 \le k \le K-1$ such that $b_k \le (1-\delta_T ) a_{k+1}$.
\end{prop}
We have thus reduced the problem of estimating $P_N(\alpha)$ to $U_N$, and the rest of the section is devoted to studying the latter sequence. Our main strategy will be to start with an arbitrary non-negative integer $N=\sum_{k=0}^{K-1} b_k q_k$, and to change its Ostrowski coefficients one by one; we call such a transformation a \textit{projection}. After finitely many projections we will transform all Ostrowski coefficients $b_k$ with $k_0 \le k \le K-1$ to $b_k^*=\lfloor (5/6) a_{k+1} \rfloor$. Keeping track of the effect of each projection, we will be able to compare $U_N$ to $U_{N^*}$.

\begin{proof}[Proof of Proposition \ref{PNUNprop}] Let $E_k(N)$ be as in Proposition \ref{pqkapproximation} if $b_k \ge 1$, and $E_k(N)=0$ if $b_k=0$. By Definition \ref{UNdef}, for any $k_0 \le k \le K-1$ we have
\[ \prod_{b=0}^{b_k-1} P_{q_k} (\alpha, (-1)^k (b q_k \| q_k \alpha \| + \varepsilon_k (N))/q_k) = u_k(N) e^{E_k(N)} , \]
and hence from the factorization \eqref{pnproductform} we get
\begin{equation}\label{PNalphaUNestimate}
P_N(\alpha ) = \left( \prod_{k=0}^{k_0-1} \prod_{b=0}^{b_k-1} P_{q_k} (\alpha, (-1)^k (b q_k \| q_k \alpha \| + \varepsilon_k (N))/q_k) \right) U_N \prod_{k=k_0}^{K-1} e^{E_k(N)} .
\end{equation}
We start by finding upper and lower bounds for the first factor independent of $N$. For an upper bound, simply use $P_{q_k} (\alpha , x) \le 2^{q_k}$ to get
\[ \prod_{k=0}^{k_0-1} \prod_{b=0}^{b_k-1} P_{q_k} (\alpha, (-1)^k (b q_k \| q_k \alpha \| + \varepsilon_k (N))/q_k) \le 2^{q_1 + \cdots +q_{k_0}} \ll_{\alpha} 1. \]
To see a lower bound, let $0 \le k \le k_0-1$ and $0 \le b \le b_k-1$. Then
\[ P_{q_k}(\alpha , (-1)^k (b q_k \| q_k \alpha \| + \varepsilon_k (N))/q_k ) = \prod_{n=1}^{q_k} |2 \sin (\pi ((n+bq_k) \alpha + (-1)^k \varepsilon_k (N) /q_k)|. \]
Here $(n+bq_k) \le q_k+(a_{k+1}-1)q_k < q_{k+1}$, and thus by the best rational approximation property and Lemma \ref{epsilonklemma},
\[ \| (n+bq_k) \alpha + (-1)^k \varepsilon_k (N)/q_k \| \ge \| (n+bq_k) \alpha \| - |\varepsilon_k (N)|/q_k \ge \| q_k \alpha \| - \| q_{k+1} \alpha \| . \]
It follows that $\| (n+bq_k) \alpha + (-1)^k \varepsilon_k (N) /q_k \| \gg_{\alpha} 1$, and hence
\[ P_{q_k}(\alpha , (-1)^k (b q_k \| q_k \alpha \| + \varepsilon_k (N))/q_k ) \gg_{\alpha} 1. \]
The first factor in \eqref{PNalphaUNestimate} is thus both $\ll_{\alpha} 1$ and $\gg_{\alpha} 1$, therefore
\[ \log P_N(\alpha ) = \log U_N + \sum_{k=k_0}^{K-1} E_k(N) + O_{\alpha} (1) . \]
By Proposition \ref{pqkapproximation} here $E_k(N) \le C/(a_{k+1}q_k)$ for all $k_0 \le k \le K-1$, and also $E_k(N) \ge -C\frac{\log (2/\delta_T)}{\delta_T^2}(1/a_{k+1} + 1/q_k^2)$ for all $k_0 \le k \le K-1$ such that $b_k \le (1-\delta_T ) a_{k+1}$ with a universal constant $C>0$. Let $F_k(N) = \max \{ -E_k(N), 0\}$ if $b_k>(1-\delta_T ) a_{k+1}$, and $F_k(N)=0$ otherwise. Then
\[ F_k (N) = -E_k(N) + O_T \left( \frac{1}{a_{k+1}} + \frac{1}{q_k^2}  \right) \]
for all $k_0 \le k \le K-1$, and the claim follows.
\end{proof}

\subsection{Key estimate for projections}

We now introduce the main technical tool in the proof of Theorem \ref{maxpntheorem}, and establish its key property.
\begin{definition}\label{projectiondef} Let $N=\sum_{k=0}^{K-1} b_k q_k$ be the Ostrowski expansion of a non-negative integer, and let $k_0 \le m \le K-1$ and $0 \le B \le (1-\delta_T )a_{m+1}$ be integers. The projection of $N$ with respect to the index $m$ and the integer $B$ is $\mathrm{proj}_{m,B} (N):=N'= \sum_{k=0}^{K-1} b_k'q_k$, where $b_k'=b_k$ for all $k \neq m$, and $b_m'=B$.
\end{definition}

\begin{prop}\label{projectionprop} Let $N=\sum_{k=0}^{K-1} b_k q_k$ be the Ostrowski expansion of a non-negative integer, and let $k_0 \le m \le K-1$ and $0 \le B \le (1-\delta_T)a_{m+1}$ be integers. Assume that $b_k \le (1-\delta_T ) a_{k+1}$ for all $k_0 \le k < m$. If $m>k_0$, then $\mathrm{proj}_{m,B} (N)=N'$ satisfies
\begin{equation}\label{logUN'-logUNinequality}
\log U_{N'} - \log U_N \ge \log u_m (N') - \log u_m (N) - (\log (b_{m-1}+1)) \frac{b_m'-b_m}{a_{m+1}} - O_T \left( \frac{|b_m'-b_m|}{a_{m+1}} \right) .
\end{equation}
If $m>k_0$ and $b_m \le (1-\delta_T ) a_{m+1}$, then \eqref{logUN'-logUNinequality} holds with equality. If $m=k_0$, then \eqref{logUN'-logUNinequality} with the term $- \log (b_{m-1}+1)\frac{b_m'-b_m}{a_{m-1}}$ removed holds with equality.
\end{prop}

\begin{proof} For the sake of readability, let $f(x)=|2 \sin (\pi x)|$. By definition \eqref{varepsilonkdef}, we have
\begin{equation}\label{varepsilonkN'-varepsilonkN}
\varepsilon_k (N')-\varepsilon_k (N) = (-1)^{k+m} (b_m'-b_m) q_k \| q_m \alpha \| \qquad \textrm{for all } k_0 \le k <m,
\end{equation}
and $\varepsilon_k (N') = \varepsilon_k (N)$ for all $m \le k \le K-1$. Recalling Definition \ref{UNdef}, it follows that $u_k(N')=u_k(N)$ for all $m<k \le K-1$, and hence
\begin{equation}\label{logun'-logun}
\log U_{N'} - \log U_N = \log u_m(N') - \log u_m (N) + \sum_{k=k_0}^{m-1} \left( \log u_k (N') - \log u_k (N) \right) .
\end{equation}
Let $k_0 \le k \le m-1$, and consider the corresponding term in the sum on the right hand side of \eqref{logun'-logun}. If $b_k'=b_k=0$, then $\log u_k(N')=\log u_k(N)=0$. Otherwise,
\begin{equation}\label{logukN'-logukN}
\begin{split} \log u_k (N') - \log u_k(N) = &\sum_{b=1}^{b_k-1} \log \frac{f(b q_k \| q_k \alpha \| + \varepsilon_k (N'))}{f(b q_k \| q_k \alpha \| + \varepsilon_k (N))} \\ &+ \sum_{b=0}^{b_k-1} \left( V_k(b q_k \| q_k \alpha \| + \varepsilon_k(N') ) - V_k(b q_k \| q_k \alpha \| + \varepsilon_k(N) ) \right) \\ &+\log \frac{b_k q_k \| q_k \alpha \| + \varepsilon_k (N')}{b_k q_k \| q_k \alpha \| + \varepsilon_k (N)} .  \end{split}
\end{equation}
To estimate the second term in \eqref{logukN'-logukN}, note that by Lemma \ref{epsilonklemma} we have
\[ b q_k \| q_k \alpha \| + \varepsilon_k (N) \ge - (1-\delta_T ), \]
and that by the assumption $b_k \le (1-\delta_T ) a_{k+1}$,
\[ b q_k \| q_k \alpha \| + \varepsilon_k (N) \le b_k q_k \| q_k \alpha \| \le 1-\delta_T . \]
Lemma \ref{Vkxlemma} (i) implies that $|V_k'(x)| \ll_T 1/a_{k+1}$ on the interval $[-(1-\delta_T ), 1-\delta_T ]$, therefore
\[ \left| \sum_{b=0}^{b_k-1} \left( V_k(b q_k \| q_k \alpha \| + \varepsilon_k(N') ) - V_k(b q_k \| q_k \alpha \| + \varepsilon_k(N) ) \right) \right| \ll_T |\varepsilon_k (N') - \varepsilon_k (N)| . \]
Using \eqref{varepsilonkN'-varepsilonkN},
\[ \sum_{k=k_0}^{m-1} |\varepsilon_k (N') - \varepsilon_k(N)| \le |b_m'-b_m| \cdot \| q_m \alpha \| \sum_{k=k_0}^{m-1} q_k \ll |b_m'-b_m| \cdot \| q_m \alpha \| q_{m-1} \le \frac{|b_m'-b_m|}{a_{m+1}} ,  \]
consequently from \eqref{logun'-logun} and \eqref{logukN'-logukN} we get
\begin{equation}\label{logun'-logun2}
\begin{split} & \log U_{N'} - \log U_N \\ = &\log u_m(N') - \log u_m (N) \\ &+\sum_{\substack{k=k_0 \\ b_k \ge 1}}^{m-1} \left( \sum_{b=1}^{b_k-1} \log \frac{f(b q_k \| q_k \alpha \| + \varepsilon_k (N'))}{f(b q_k \| q_k \alpha \| + \varepsilon_k (N))} + \log \frac{b_k q_k \| q_k \alpha \| + \varepsilon_k (N')}{b_k q_k \| q_k \alpha \| + \varepsilon_k (N)}  \right) \\ &+ O_T \left( \frac{|b_m'-b_m|}{a_{m+1}} \right) . \end{split}
\end{equation}
Next, we show that the sum over $k_0 \le k \le m-2$ in the previous formula is negligible. Let $k_0 \le k \le m-2$. By assumption, we have $b_{k+1} \le (1-\delta_T ) a_{k+2}$, and so $\varepsilon_k (N) \ge -(1-\delta_T/3 ) q_k \| q_k \alpha \|$ and $\varepsilon_k (N') \ge -(1-\delta_T/3 ) q_k \| q_k \alpha \|$ follow from Lemma \ref{epsilonklemma}. In particular, $bq_k \| q_k \alpha \| +\varepsilon_k (N)$ and $bq_k \| q_k \alpha \| +\varepsilon_k (N')$ both lie in the interval $[ (b-1+\delta_T/3 ) q_k \| q_k \alpha \| , 1-\delta_T ]$ for all $1 \le b \le b_k-1$. Since $|(\log f(x))'| \ll_T \frac{1}{(b-1+\delta_T /3)q_k \| q_k \alpha \|}$ on this interval, we have
\[ \begin{split} \sum_{b=1}^{b_k-1} \left| \log \frac{f(b q_k \| q_k \alpha \| + \varepsilon_k (N'))}{f(b q_k \| q_k \alpha \| + \varepsilon_k (N))} \right| &\ll_T \sum_{b=1}^{b_k-1} \frac{|\varepsilon_k (N') - \varepsilon_k(N)|}{(b-1+\delta_T /3) q_k \| q_k \alpha \|} \\ &\ll_T a_{k+1} (\log a_{k+1}) |\varepsilon_k (N') - \varepsilon_k (N)| \\ &\ll_T a_{k+1} a_{k+2} q_k |b_m'-b_m| \cdot \| q_m \alpha \| \\ &\le q_{k+2} |b_m'-b_m| \cdot \| q_m \alpha \| , \end{split} \]
where we used \eqref{varepsilonkN'-varepsilonkN} and condition \eqref{logak/ak+1}. Since $b_k q_k \| q_k \alpha \| + \varepsilon_k(N)$ and $b_k q_k \| q_k \alpha \| + \varepsilon_k(N')$ both lie in the interval $[(\delta_T/3 ) q_k \| q_k \alpha \|, 2]$, and $|(\log x)'| \ll_T a_{k+1}$ on this interval, we similarly get
\[ \left| \log \frac{b_k q_k \| q_k \alpha \| + \varepsilon_k (N')}{b_k q_k \| q_k \alpha \| + \varepsilon_k (N)} \right| \ll_T a_{k+1}|\varepsilon_k(N') - \varepsilon_k (N)| \ll q_{k+1} |b_m' - b_m| \cdot \| q_m \alpha \| . \]
From the previous two formulas and $\sum_{k=k_0}^{m-2} q_{k+2} \ll q_m$ we get
\[ \sum_{\substack{k=k_0 \\ b_k \ge 1}}^{m-2} \left| \sum_{b=1}^{b_k-1} \log \frac{f(b q_k \| q_k \alpha \| + \varepsilon_k (N'))}{f(b q_k \| q_k \alpha \| + \varepsilon_k (N))} + \log \frac{b_k q_k \| q_k \alpha \| + \varepsilon_k (N')}{b_k q_k \| q_k \alpha \| + \varepsilon_k (N)}  \right| \ll_T \frac{|b_m'-b_m|}{a_{m+1}} , \]
and hence if $m>k_0$, \eqref{logun'-logun2} simplifies as
\begin{equation}\label{logun'-logun3}
\begin{split} \log U_{N'} - \log U_N = &\log u_m (N') - \log u_m (N) \\ &+ \sum_{b=1}^{b_{m-1}-1} \log \frac{f(b q_{m-1} \| q_{m-1} \alpha \| + \varepsilon_{m-1} (N'))}{f(b q_{m-1} \| q_{m-1} \alpha \| + \varepsilon_{m-1} (N))} \\ &+I_{\{b_{m-1} \ge 1\}} \log \frac{b_{m-1} q_{m-1} \| q_{m-1} \alpha \| + \varepsilon_{m-1} (N')}{b_{m-1} q_{m-1} \| q_{m-1} \alpha \| + \varepsilon_{m-1} (N)} \\ & + O_T \left( \frac{|b_m'-b_m|}{a_{m+1}} \right) . \end{split}
\end{equation}
If $m=k_0$, then \eqref{logun'-logun3} holds with the second and third terms on the right hand side removed, and the claim for $m=k_0$ follows.

Let $m>k_0$. To proceed, we distinguish between two cases: Case 1 is $\varepsilon_{m-1}(N)<-(1-\delta_T /3) q_m \| q_m \alpha \|$, and Case 2 is $\varepsilon_{m-1}(N) \ge -(1-\delta_T /3) q_m \| q_m \alpha \|$. We will show that \eqref{logUN'-logUNinequality} holds in Case 1, and that \eqref{logUN'-logUNinequality} holds with equality in Case 2. Note that this will prove the proposition; indeed, \eqref{logUN'-logUNinequality} follows in either case, whereas by Lemma \ref{epsilonklemma} the additional assumption $b_m \le (1-\delta_T ) a_{m+1}$ ensures that we are in Case 2.

\vspace{5mm}

\noindent\textbf{Case 1.} Assume that $\varepsilon_{m-1}(N)<-(1-\delta_T /3) q_m \| q_m \alpha \|$. By assumption, $b_m' \le (1-\delta_T) a_{m+1}$, and hence $\varepsilon_{m-1} (N') \ge -(1-\delta_T /3) q_m \| q_m \alpha \|$ follows from Lemma \ref{epsilonklemma}. In particular, $\varepsilon_{m-1}(N) < \varepsilon_{m-1}(N')$, and hence
\begin{equation}\label{Ibm-1lowerbound}
I_{\{b_{m-1} \ge 1\}} \log \frac{b_{m-1} q_{m-1} \| q_{m-1} \alpha \| + \varepsilon_{m-1} (N')}{b_{m-1} q_{m-1} \| q_{m-1} \alpha \| + \varepsilon_{m-1} (N)} \ge 0.
\end{equation}
If $b_{m-1}=0$ or $1$, then \eqref{logUN'-logUNinequality} follows from \eqref{logun'-logun3}; therefore we may assume that $b_{m-1} \ge 2$. Since $(\log f(x))'=\pi \cot (\pi x)$, for any $1 \le b \le b_{m-1}-1$ we have
\[ \begin{split} \log &\frac{f(b q_{m-1} \| q_{m-1} \alpha \| + \varepsilon_{m-1} (N'))}{f(b q_{m-1} \| q_{m-1} \alpha \| + \varepsilon_{m-1} (N))} \\ & \hspace{25mm} = \int_{b q_{m-1} \| q_{m-1} \alpha \| + \varepsilon_{m-1} (N)}^{b q_{m-1} \| q_{m-1} \alpha \| + \varepsilon_{m-1} (N')} \pi \cot (\pi x) \, \mathrm{d} x \\ &\hspace{25mm} \ge \pi \cot (\pi (b q_{m-1} \| q_{m-1} \alpha \| + \varepsilon_{m-1} (N'))) \cdot (\varepsilon_{m-1} (N') - \varepsilon_{m-1}(N)) . \end{split} \]
Using $|\pi \cot (\pi x) -1/x| \ll_T 1$ for all $|x| \le 1-\delta_T$, we get
\[ \begin{split} & \sum_{b=1}^{b_{m-1}-1} \pi \cot (\pi (b q_{m-1} \| q_{m-1} \alpha \| + \varepsilon_{m-1} (N'))) \\ = & \sum_{b=1}^{b_{m-1}-1} \frac{1}{b q_{m-1} \| q_{m-1} \alpha \| + \varepsilon_{m-1}(N')} + O_T(a_m) \\ = & a_m \log (b_{m-1}+1) + O_T(a_m) . \end{split} \]
The previous two formulas give
\[ \begin{split} \sum_{b=1}^{b_{m-1}-1}  \log &\frac{f(b q_{m-1} \| q_{m-1} \alpha \| + \varepsilon_{m-1} (N'))}{f(b q_{m-1} \| q_{m-1} \alpha \| + \varepsilon_{m-1} (N))} \\ &\ge a_m (\log (b_{m-1}+1)) (\varepsilon_{m-1} (N') - \varepsilon_{m-1}(N)) + O_T \left( a_m |\varepsilon_{m-1} (N') - \varepsilon_{m-1}(N)| \right) \\ &= - (\log (b_{m-1}+1)) \frac{b_m'-b_m}{a_{m+1}} + O_T \left( \frac{|b_m'-b_m|}{a_{m+1}} \right) , \end{split} \]
therefore \eqref{logun'-logun3} and \eqref{Ibm-1lowerbound} imply the desired inequality \eqref{logUN'-logUNinequality}. This concludes the proof of Case 1.

\vspace{5mm}

\noindent\textbf{Case 2.} Assume that $\varepsilon_{m-1} (N) \ge -(1-\delta_T /3) q_m \| q_m \alpha \|$. Repeating arguments from above, we now have
\begin{equation}\label{Ibm-1bound}
I_{\{b_{m-1} \ge 1\}} \log \frac{b_{m-1} q_{m-1} \| q_{m-1} \alpha \| + \varepsilon_{m-1} (N')}{b_{m-1} q_{m-1} \| q_{m-1} \alpha \| + \varepsilon_{m-1} (N)} = O_T \left( \frac{|b_m'-b_m|}{a_{m+1}} \right) .
\end{equation}
If $b_{m-1}=0$ or $1$, then the sum in \eqref{logun'-logun3} is empty, and it follows that \eqref{logUN'-logUNinequality} holds with equality; thus we may again assume that $b_{m-1} \ge 2$. Since we are in Case 2, for any $1 \le b \le b_{m-1}-1$ the points $b q_{m-1} \| q_{m-1} \alpha \| + \varepsilon_{m-1}(N)$ and $b q_{m-1} \| q_{m-1} \alpha \| + \varepsilon_{m-1}(N')$ lie in the interval
\[ [(b-1+\delta_T /3)q_{m-1} \| q_{m-1} \alpha \|, (b+1) q_{m-1} \| q_{m-1} \alpha \| ] . \]
Note that here $(b+1)q_{m-1} \| q_{m-1} \alpha \| \le 1-\delta_T$. We have $(\log f(x))'=\pi \cot (\pi x)$, and $|(\log f(x))''| = |\pi/\sin^2 (\pi x)| \ll_T \frac{1}{b^2 q_{m-1}^2 \| q_{m-1} \alpha \|^2}$ on the same interval. Applying a second order Taylor Formula, we thus get
\[ \begin{split} \log &\frac{f(b q_{m-1} \| q_{m-1} \alpha \| + \varepsilon_{m-1}(N'))}{f(b q_{m-1} \| q_{m-1} \alpha \| + \varepsilon_{m-1}(N) )} \\ & \hspace{25mm} = \pi \cot (\pi (b q_{m-1} \| q_{m-1} \alpha \| + \varepsilon_{m-1}(N'))) (\varepsilon_{m-1}(N') - \varepsilon_{m-1}(N)) \\ & \hspace{25mm} + O_T \left( \frac{(\varepsilon_{m-1}(N') - \varepsilon_{m-1}(N))^2}{b^2 q_{m-1}^2 \| q_{m-1} \alpha \|^2} \right) . \end{split} \]
The contribution of the error term is negligible:
\[ \begin{split} \sum_{b=1}^{b_{m-1}-1} \frac{(\varepsilon_{m-1}(N') - \varepsilon_{m-1}(N))^2}{b^2 q_{m-1}^2 \| q_{m-1} \alpha \|^2} &\ll a_m^2 (\varepsilon_{m-1}(N') - \varepsilon_{m-1}(N))^2 \\ &= a_m^2 q_{m-1}^2 (b_m'-b_m)^2 \| q_m \alpha \|^2 \le \frac{|b_m'-b_m|}{a_{m+1}} . \end{split} \]
Arguing as in Case 1, the previous two formulas yield
\[ \sum_{b=1}^{b_{m-1}-1}  \log \frac{f(b q_{m-1} \| q_{m-1} \alpha \| + \varepsilon_{m-1} (N'))}{f(b q_{m-1} \| q_{m-1} \alpha \| + \varepsilon_{m-1} (N))} = - (\log (b_{m-1}+1)) \frac{b_m'-b_m}{a_{m+1}} + O_T \left( \frac{|b_m'-b_m|}{a_{m+1}} \right) , \]
therefore \eqref{logun'-logun3} and \eqref{Ibm-1bound} imply that \eqref{logUN'-logUNinequality} holds with equality. This concludes the proof of Case 2.
\end{proof}

\subsection{Regularizing and optimizing projections}\label{regoptsection}

We introduced the concept of a projection in Definition \ref{projectiondef}. Starting with a non-negative integer with Ostrowski expansion $N=\sum_{k=0}^{K-1}b_kq_k$, our strategy is to apply projections to $N$ in two rounds. In the first round we project the coefficients with $b_k>(1-\delta_T )a_{k+1}$ to
\begin{equation}\label{bk**def}
b_k^{**}:= \left\{ \begin{array}{ll} 0 & \textrm{if } a_{k+1}=2, \\ \lfloor (1-\delta_T ) a_{k+1} \rfloor & \textrm{if } a_{k+1} \neq 2 \end{array} \right.
\end{equation}
in increasing order of the indices $k_0 \le k \le K-1$. We call such a transformation a \textit{regularizing projection}; its aim is to get away from the singularity of $\log |2 \sin (\pi x)|$ at $x=1$. We note that the special value of $b_k^{**}$ in the case $a_{k+1}=2$ ($0$ instead of $1$) serves a technical purpose, and will not cause difficulties in the end. After the first round of projections $N$ is transformed into an integer whose $k$-th Ostrowski coefficient is $\le (1-\delta_T )a_{k+1}$ for all $k_0 \le k \le K-1$. As we will see, the value of $U_N$ does not decrease up to a small error during the first round.

In the second round we project each Ostrowski coefficient to $b_k^*=\lfloor (5/6) a_{k+1} \rfloor$, in increasing order of the indices $k_0 \le k \le K-1$. We call such a transformation an \textit{optimizing projection}. We now estimate the effect of a projection on the value of $U_N$ based on its type.
\begin{prop}\label{regoptprojections} Let $N=\sum_{k=0}^{K-1} b_k q_k$ be the Ostrowski expansion of a non-negative integer, and let $k_0 \le m \le K-1$.
\begin{enumerate}
\item[(i)] (Regularizing projection.) Assume that $b_k \le (1-\delta_T ) a_{k+1}$ for all $k_0 \le k <m$, and that $b_m>(1-\delta_T ) a_{m+1}$. Then $\mathrm{proj}_{m,b_m^{**}}(N)=N'$ satisfies
\[ \log U_{N'} - \log U_N \ge -O_T(1/a_{m+1}) -O_{\alpha} (1/q_m) . \]
\item[(ii)] (Optimizing projection.) Assume that $b_k \le (1-\delta_T ) a_{k+1}$ for all $k_0 \le k \le K-1$, and that $b_k=b_k^*$ for all $k_0 \le k <m$. Then $\mathrm{proj}_{m,b_m^*}(N)=N'$ satisfies
\begin{equation}\label{optprojapproximation}
\begin{split} \log U_{N'} - \log U_N = &a_{m+1} \int_{b_m/a_{m+1}}^{b_m'/a_{m+1}} \log |2 \sin (\pi x)| \, \mathrm{d}x \\ &+ O_T \left( \frac{|b_m'-b_m|}{a_{m+1}} +I_{\{ b_m \le 0.01 a_{m+1} \}} \log a_{m+1} \right) +O_{\alpha}(1/q_m) , \end{split}
\end{equation}
and also
\begin{equation}\label{optprojlowerbound}
\log U_{N'} - \log U_N \ge 0.2326 \frac{(b_m'-b_m)^2}{a_{m+1}} -O_T \left( 1/a_{m+1} \right) -O_{\alpha} (1/q_m) .
\end{equation}
\end{enumerate}
\end{prop}

We first prove a lemma which will help in the case $a_{m+1} \ll 1$, and then give the proofs of Proposition \ref{regoptprojections} (i) and (ii).

\begin{lem}\label{Alemma} Let $N=\sum_{k=0}^{K-1} b_k q_k$ be the Ostrowski expansion of a non-negative integer, and let $k_0 \le m \le K-1$. If $a_{m+1} \le A$ with some constant $A \ge 1$, then
\begin{equation} \label{with_equ}
\log u_m(N) \le O_{T,A}(1) + O_{\alpha} (1/q_m) . 
\end{equation}
If in addition $b_m \le (1-\delta_T ) a_{m+1}$ and $b_{m+1} \le (1-\delta_T ) a_{m+2}$, then \eqref{with_equ} holds with $|\log u_m(N)|$ instead of $\log u_m(N)$ on the left-hand side.
\end{lem}

\begin{proof} Recall Definition \ref{UNdef}. If $b_m=0$, then $u_m(N)=1$, and we are done. We may thus assume that $b_m \ge 1$.

First, we prove the upper bound. From the definition of $u_m(N)$ we immediately see that
\[ u_m(N) \le 2^A \exp \left( \sum_{b=0}^{b_m-1} V_m (b q_m \| q_m \alpha \| + \varepsilon_m(N)) \right) 4 \pi .  \]
Lemma \ref{Vkxlemma} (i) implies that $V_m$ is decreasing on $(-1,1)$, and $|V_m'(x)|\ll_T 1/a_{m+1}$ on $[-(1-\delta_T ),0]$. It is then readily seen that
\[ \sum_{b=0}^{b_m-1} V_m (b q_m \| q_m \alpha \| + \varepsilon_m (N)) \le b_m V_m(0) + O_T(1) = O_{T,A}(1) + O_{\alpha} (1/q_m) , \]
and the claim $\log u_m(N) \le O_{T,A}(1) + O_{\alpha} (1/q_m)$ follows.

Next, assume in addition, that $b_m \le (1-\delta_T ) a_{m+1}$ and $b_{m+1} \le (1-\delta_T ) a_{m+2}$. By Lemma \ref{epsilonklemma} we then have $\varepsilon_m (N) \ge -(1-\delta_T /3) q_m \| q_m \alpha \|$. Therefore the points $bq_m \| q_m \alpha \| + \varepsilon_m(N)$, $1 \le b \le b_m-1$ are bounded away from $0$ and $1$, and hence
\[ \prod_{b=1}^{b_m-1} |2 \sin (\pi (b q_m \| q_m \alpha \| + \varepsilon_m (N)))| \gg_{T,A} 1 . \]
Similarly, since the points $b q_m \| q_m \alpha \| + \varepsilon_m (N)$, $0 \le b \le b_m-1$ lie in the interval $[-(1-\delta_T ), 1-\delta_T ]$ and $|V_m'(x)| \ll_T 1/a_{m+1}$ on this interval by Lemma \ref{Vkxlemma} (i), we have
\[ \left| \sum_{b=0}^{b_m-1} V_m (b q_m \| q_m \alpha \| + \varepsilon_m (N)) \right| = b_m |V_m(0)| + O_T(1) = O_{T,A}(1) + O_{\alpha} (1/q_m) . \]
Finally, note that
\[ 2 \pi (b_m q_m \|q_m \alpha \| + \varepsilon_m (N)) \ge 2 \pi \cdot (\delta_T /3) q_m \| q_m \alpha \| \gg_{T,A} 1. \]
The previous three estimates and the definition of $u_m(N)$ show that $\log u_m(N) \ge -O_{T,A}(1)-O_{\alpha} (1/q_m)$, as claimed.
\end{proof}

\begin{proof}[Proof of Proposition \ref{regoptprojections} (i)] For the sake of readability, let $f(x)=|2 \sin (\pi x)|$ and $\varepsilon_m = \varepsilon_m(N)=\varepsilon_m(N')$. Since $b_m>b_m'$, and in particular $-\log (b_{m-1}+1) \frac{b_m'-b_m}{a_{m+1}} \ge 0$, Proposition \ref{projectionprop} gives
\begin{equation}\label{logUN'-logUN(i)}
\log U_{N'} - \log U_N \ge \log u_m(N') - \log u_m(N) -O_T \left( \frac{|b_m'-b_m|}{a_{m+1}} \right) .
\end{equation}
Let $A \ge 1$ be a large constant depending only on $T$ (and $\delta_T$), to be chosen. We will distinguish between three cases depending on the size of $a_{m+1}$.

\vspace{5mm}

\noindent\textbf{Case 1.} Assume that $a_{m+1}=1$ or $2$. Then by construction $b_m'=b_m^{**}=0$, and hence $u_m(N')=1$. From \eqref{logUN'-logUN(i)} and Lemma \ref{Alemma} we thus get $\log U_{N'} - \log U_N \ge - O_T(1) - O_{\alpha} (1/q_m)$, and the claim follows.

\vspace{5mm}

\noindent\textbf{Case 2.} Assume that $3 \le a_{m+1} \le A$. Then by construction $b_m'=b_m^{**} \ge 2$. Recalling Definition \ref{UNdef}, we have
\[ \begin{split} \log u_m(N') - \log u_m (N) = &- \sum_{b=b_m'}^{b_m-1} \log f (b q_m \| q_m \alpha \| + \varepsilon_m ) - \sum_{b=b_m'}^{b_m-1} V_m (b q_m \| q_m \alpha \| + \varepsilon_m ) \\ &+\log \frac{b_m' q_m \| q_m \alpha \| + \varepsilon_m}{b_m q_m \| q_m \alpha \| + \varepsilon_m} . \end{split} \]
Since $-\log f(x)$ is bounded from below, the first term is $\ge -O_A (1)$. Similar to the proof of Lemma \ref{Alemma}, it is easy to see that the second term is $\ge -O_{T,A}(1)-O_{\alpha}(1/q_m)$. Finally, the last term is
\[ \log \frac{b_m' q_m \| q_m \alpha \| + \varepsilon_m}{b_m q_m \| q_m \alpha \| + \varepsilon_m} \ge \log \frac{b_m'-1}{b_m+1} \ge -O_T(1). \]
Hence $\log u_m(N') - \log u_m(N) \ge -O_{T,A}(1) - O_{\alpha} (1/q_m)$, and the claim follows from \eqref{logUN'-logUN(i)}.

\vspace{5mm}

\noindent\textbf{Case 3.} Assume that $a_{m+1}>A$. By Definition \ref{UNdef}, we again have
\[ \begin{split} \log u_m(N') - \log u_m (N) = &- \sum_{b=b_m'}^{b_m-1} \log f (b q_m \| q_m \alpha \| + \varepsilon_m ) - \sum_{b=b_m'}^{b_m-1} V_m (b q_m \| q_m \alpha \| + \varepsilon_m ) \\ &+\log \frac{b_m' q_m \| q_m \alpha \| + \varepsilon_m}{b_m q_m \| q_m \alpha \| + \varepsilon_m} . \end{split} \]
From Lemma \ref{Vkxlemma}, in particular from the fact that $V_m(x)$ is decreasing, we get
\[ \begin{split} - \sum_{b=b_m'}^{b_m-1} V_m (b q_m \| q_m \alpha \| + \varepsilon_m ) &\ge - (b_m-b_m') V_m(0) \\ &= - (b_m-b_m') \frac{\log a_m}{a_{m+1}} - O_T \left( \frac{|b_m'-b_m|}{a_{m+1}} \right) - O_{\alpha} \left( \frac{|b_m'-b_m|}{q_{m+1}} \right) \\ &\ge - T(b_m-b_m') - O_T \left( \frac{|b_m'-b_m|}{a_{m+1}} \right) - O_{\alpha} (1/q_m) . \end{split} \]
Since $b_m q_m \| q_m \alpha \| + \varepsilon_k$ and $b_m' q_m \| q_m \alpha \| + \varepsilon_k$ are bounded away from zero, we also have
\[ \log \frac{b_m' q_m \| q_m \alpha \| + \varepsilon_m}{b_m q_m \| q_m \alpha \| + \varepsilon_m} \ll \frac{|b_m-b_m'|}{a_{m+1}} . \]
By the previous two formulas, \eqref{logUN'-logUN(i)} simplifies to
\begin{equation}\label{logUN'-logUN(i)2}
\log U_{N'} - \log U_N \ge - \sum_{b=b_m'}^{b_m-1} \log f (b q_m \| q_m \alpha \| + \varepsilon_m ) - T(b_m-b_m') -O_T \left( \frac{|b_m'-b_m|}{a_{m+1}} \right) - O_{\alpha} (1/q_m) .
\end{equation}
For all $b_m' \le b \le b_m-1$,
\[ b q_m \| q_m \alpha \| + \varepsilon_k (N) \ge (\lfloor (1-\delta_T ) a_{m+1} \rfloor -1) q_m \| q_m \alpha \| \ge 1-2\delta_T \]
provided that $A$ is large enough in terms of $T$ and $\delta_T$. Choosing $\delta_T \le 1/(4 \pi e^{2T})$ ensures that $\log f(x) \le -2T$ on the interval $[1-2 \delta_T, 1)$. Hence every term in the sum in \eqref{logUN'-logUN(i)2} is $\le -2T$, and we get
\[ \log U_{N'} - \log U_N \ge T (b_m-b_m') -O_T \left( \frac{|b_m'-b_m|}{a_{m+1}} \right) - O_{\alpha} (1/q_m) . \]
Choosing $A$ large enough in terms of $T$ and $\delta_T$, the second error term is negligible compared to $T(b_m-b_m')$. Hence $\log U_{N'} - \log U_N \ge -O_{\alpha}(1/q_m)$, and the claim follows.
\end{proof}

\begin{proof}[Proof of Proposition \ref{regoptprojections} (ii)] Again, let $f(x)=|2 \sin (\pi x)|$ and $\varepsilon_m = \varepsilon_m(N)=\varepsilon_m(N')$. If $m=k_0$, then both sides of \eqref{optprojapproximation} and \eqref{optprojlowerbound} are $O_{\alpha} (1)$, and we are done. We may thus assume that $m>k_0$. By the assumption $b_{m-1}=b_{m-1}^*=\lfloor (5/6)a_m \rfloor$ we have
\[ - (\log (b_{m-1}+1)) \frac{b_m'-b_m}{a_{m+1}} = -(b_m'-b_m) \frac{\log a_m}{a_{m+1}} + O \left( \frac{|b_m'-b_m|}{a_{m+1}} \right) , \]
hence Proposition \ref{projectionprop} now gives
\begin{equation}\label{logUN'-logUN(ii)}
\log U_{N'} - \log U_N = \log u_m (N') - \log u_m (N) - (b_m'-b_m) \frac{\log a_m}{a_{m+1}} - O_T \left( \frac{|b_m'-b_m|}{a_{m+1}} \right) .
\end{equation}
Let $A \ge 1$ be a large constant depending only on $T$, to be chosen. We distinguish between four cases.

\vspace{5mm}

\noindent\textbf{Case 1.} Assume that $a_{m+1} \le A$. By Lemma \ref{Alemma}, both $\log u_m(N')$ and $\log u_m(N)$ are $O_{T,A}(1)+O_{\alpha}(1/q_m)$. Since $|b_m'-b_m| \frac{\log a_m}{a_{m+1}} \ll_{T,A} 1$, from \eqref{logUN'-logUN(ii)} we get
\[ \log U_{N'} - \log U_N = O_{T,A}(1) + O_{\alpha} (1/q_m),  \]
and the claims \eqref{optprojapproximation} and \eqref{optprojlowerbound} follow.

\vspace{5mm}

\noindent\textbf{Case 2.} Assume that $a_{m+1}>A$ and $1 \le b_m \le b_m'$. From Definition \ref{UNdef} we now get
\[ \begin{split} \log u_m(N') - \log u_m (N) = &\sum_{b=b_m}^{b_m'-1} \log f (b q_m \| q_m \alpha \| + \varepsilon_m ) + \sum_{b=b_m}^{b_m'-1} V_m (b q_m \| q_m \alpha \| + \varepsilon_m ) \\ &+\log \frac{b_m' q_m \| q_m \alpha \| + \varepsilon_m}{b_m q_m \| q_m \alpha \| + \varepsilon_m} . \end{split} \]
The assumption $b_k \le (1-\delta_T ) a_{k+1}$ for all $k_0 \le k \le K-1$ and Lemma \ref{epsilonklemma} imply that $\varepsilon_m \ge -(1-\delta_T /3) q_m \| q_m \alpha \|$. We claim that the last term satisfies
\[ \log \frac{b_m' q_m \| q_m \alpha \| + \varepsilon_m}{b_m q_m \| q_m \alpha \| + \varepsilon_m} \ll_T \frac{|b_m'-b_m|}{a_{m+1}} + I_{\{ b_m \le 0.01 a_{m+1} \}} \log a_{m+1} . \]
Indeed, if $b_m > 0.01 a_{m+1}$, then the points $b_m q_m \| q_m \alpha \| + \varepsilon_m$ and $b_m' q_m \| q_m \alpha \| + \varepsilon_m$ lie in an interval bounded away from zero, and $|(\log x)'| \ll_T 1$ on such an interval; the upper bound $|b_m'-b_m|/a_{m+1}$ follows. If $b_m \le 0.01 a_{m+1}$, then
\[ 0 \le \log \frac{b_m' q_m \| q_m \alpha \| + \varepsilon_m}{b_m q_m \| q_m \alpha \| + \varepsilon_m} \le \log \frac{2}{(\delta_T /3) q_m \| q_m \alpha \|} \ll_T \log a_{m+1}, \]
and the claimed upper bound follows once again.

Observe also that for all $b_m \le b \le b_m'-1$, the point $b q_m \| q_m \alpha \| + \varepsilon_m$ lies in $[0,5/6]$. Using Lemma \ref{Vkxlemma} (iii) we thus deduce
\[ \sum_{b=b_m}^{b_m'-1} V_m (b q_m \| q_m \alpha \| + \varepsilon_m ) = (b_m'-b_m) \frac{\log a_m}{a_{m+1}} +O_T \left( \frac{|b_m'-b_m|}{a_{m+1}} \right) + O_{\alpha} (1/q_m) , \]
hence
\[ \begin{split} \log u_m(N') - \log u_m(N) = &\sum_{b=b_m}^{b_m'-1} \log f (b q_m \| q_m \alpha \| + \varepsilon_m ) + (b_m'-b_m) \frac{\log a_m}{a_{m+1}} \\ &+ O_T \left( \frac{|b_m'-b_m|}{a_{m+1}} + I_{\{ b_m \le 0.01 a_{m+1} \}} \log a_{m+1} \right) + O_{\alpha} (1/q_m) . \end{split} \]
With a remarkable cancellation of $(b_m'-b_m) \frac{\log a_m}{a_{m+1}}$, equation \eqref{logUN'-logUN(ii)} thus simplifies to
\[ \begin{split} \log U_{N'} - \log U_N = &\sum_{b=b_m}^{b_m'-1} \log f (b q_m \| q_m \alpha \| + \varepsilon_m ) \\ &+ O_T \left( \frac{|b_m'-b_m|}{a_{m+1}} + I_{\{ b_m \le 0.01 a_{m+1} \}} \log a_{m+1} \right) + O_{\alpha} (1/q_m) . \end{split} \]

We now prove \eqref{optprojapproximation}. Assume first, that $b_m >0.01 a_{m+1}$. Then for all $b_m \le b \le b_m'-1$, the point $bq_m \| q_m \alpha \| + \varepsilon_m$ lies in an interval bounded away from $0$ and $1$. Using $q_m \| q_m \alpha \| = 1/a_{m+1}+O(1/a_{m+1}^2)$ and $|\varepsilon_m| \le 1/a_{m+1}$ we deduce
\[ |\log f(b q_m \| q_m \alpha \| + \varepsilon_m ) - \log f(b/a_{m+1}) | \ll 1/a_{m+1} , \]
and hence
\[ \sum_{b=b_m}^{b_m'-1} \log f (b q_m \| q_m \alpha \| + \varepsilon_m ) = \sum_{b=b_m}^{b_m'-1} \log f(b/a_{m+1}) +O \left( \frac{|b_m'-b_m|}{a_{m+1}} \right) . \]
Since $|\log f(x)| \ll |x-5/6|$ on our interval bounded away from 0 and 1, each term in the previous sum is $\ll |b_m-b_m'|/a_{m+1}$. Therefore by interpreting the sum as a Riemann sum,
\[ \sum_{b=b_m}^{b_m'-1} \log f(b/a_{m+1}) = a_{m+1} \int_{b_m/a_{m+1}}^{b_m'/a_{m+1}} \log f(x) \, \mathrm{d} x + O \left( \frac{|b_m'-b_m|}{a_{m+1}} \right) , \]
and \eqref{optprojapproximation} follows provided that $b_m>0.01 a_{m+1}$. If $b_m \le 0.01 a_{m+1}$, then for any $b_m \le b \le b_m'-1$, the point $bq_m \| q_m \alpha \| +\varepsilon_m$ lies in $[(b-1+\delta_T /3) q_m \| q_m \alpha \|, 5/6]$. Since $|(\log f(x))'|=|\pi \cot (\pi x)| \ll_T a_{m+1}/b$ on this interval, we now have
\[ \sum_{b=b_m}^{b_m'-1} |\log f(b q_m \| q_m \alpha \| + \varepsilon_m ) - \log f(b/a_{m+1}) | \ll_T \sum_{b=b_m}^{b_m'-1} \frac{1}{b} \ll \log a_{m+1} . \]
Note that $|\log f(b/a_{m+1})| \ll \log a_{m+1}$, hence by interpreting the sum as a Riemann sum, we now have
\[ \sum_{b=b_m}^{b_m'-1} \log f(b/a_{m+1}) = a_{m+1} \int_{b_m/a_{m+1}}^{b_m'/a_{m+1}} \log f(x) \, \mathrm{d} x + O_T \left( \log a_{m+1} \right) , \]
and \eqref{optprojapproximation} follows in the case $b_m \le 0.01 a_{m+1}$ as well.

Next, we deduce \eqref{optprojlowerbound} from \eqref{optprojapproximation}. Choosing $A$ large enough in terms of $T$, the error term in \eqref{optprojapproximation} is negligible compared to the main term (in both cases $b_m \le 0.01 a_{m+1}$ and $b_m>0.01 a_{m+1}$). Elementary calculations show that
\[ \min_{y \in [0,5/6)} \frac{1}{(5/6-y)^2} \int_y^{5/6} \log f(x) \, \mathrm{d}x = \frac{1}{(5/6)^2} \int_0^{5/6} \log f(x) \, \mathrm{d}x = 0.23260748\dots . \]
Indeed, the left hand side is an increasing function of $y$ on $[0,5/6)$; to see that its derivative is non-negative, it is enough to check that
\[ \frac{1}{5/6-y} \int_{y}^{5/6} \log f(x) \, \mathrm{d}x \ge \frac{\log f(y)}{2} , \]
and this follows from the concavity of $\log f(x)$. Therefore up to a negligible $O(1/a_{m+1})$ error in the numerical constants,
\[ a_{m+1} \int_{b_m/a_{m+1}}^{b_m'/a_{m+1}} \log f(x) \, \mathrm{d}x \ge a_{m+1} \cdot 0.2326 (b_m/a_{m+1} -5/6)^2 \ge 0.2326 \frac{(b_m'-b_m)^2}{a_{m+1}} . \]
This finishes the proof of \eqref{optprojapproximation} and \eqref{optprojlowerbound} in Case 2.

\vspace{5mm}

\noindent\textbf{Case 3.} Assume that $a_{m+1} >A$ and $b_m=0$. Then $\log u_m(N)=0$, and
\[ \begin{split} \log u_m(N')-\log u_m(N) = &\sum_{b=1}^{b_m'-1} \log f (b q_m \| q_m \alpha \| + \varepsilon_m) + \sum_{b=0}^{b_m'-1} V_m (b q_m \| q_m \alpha \| + \varepsilon_m ) \\ &+ \log (b_m' q_m \| q_m \alpha \| + \varepsilon_m ) . \end{split} \]
Following the steps of Case 2 (observe that in the first sum the summation now starts at $b=1$ instead of $b=b_m=0$), we get
\[ \log U_{N'} - \log U_N = a_{m+1} \int_{1/a_{m+1}}^{b_m'/a_{m+1}} \log f(x) \, \mathrm{d}x + O_T(\log a_{m+1}) +O_{\alpha}(1/q_m) . \]
The error of replacing the lower limit of integration by zero is negligible:
\[ a_{m+1} \int_0^{1/a_{m+1}} \log f(x) \, \mathrm{d}x \ll \log a_{m+1}, \]
and \eqref{optprojapproximation} follows. We deduce \eqref{optprojlowerbound} from \eqref{optprojapproximation} as in Case 2.

\vspace{5mm}

\noindent\textbf{Case 4.} Assume that $a_{m+1} > A$ and $b_m' \le b_m \le (1-\delta_T )a_{m+1}$. Working on the interval $[5/6,1-\delta_T]$ instead of $[0,5/6]$, the proof of \eqref{optprojapproximation} is entirely analogous to that in Case 2. Deducing \eqref{optprojlowerbound} from \eqref{optprojapproximation} is even simpler. Indeed, note that by concavity $\log f(x) \le -\pi \sqrt{3} (x-5/6)$, the right hand side being a tangent line. Hence up to a negligible $O(1/a_{m+1})$ error in the numerical constants,
\[ a_{m+1} \int_{b_m/a_{m+1}}^{b_m'/a_{m+1}} \log f(x) \, \mathrm{d}x \ge a_{m+1} \int_{b_m'/a_{m+1}}^{b_m/a_{m+1}} \pi \sqrt{3} (x-5/6) \, \mathrm{d}x \ge \frac{\pi \sqrt{3}}{2} \cdot \frac{(b_m'-b_m)^2}{a_{m+1}} . \]
The lower bound \eqref{optprojlowerbound} thus follows, in fact with the better numerical constant $\pi \sqrt{3}/2 \approx 2.72$.
\end{proof}

\subsection{Completing the proof}

\begin{proof}[Proof of Theorem \ref{maxpntheorem}] Let $N=\sum_{k=0}^{K-1} b_k q_k$ be the Ostrowski expansion of a non-negative integer, and let $N^*=\sum_{k=0}^{K-1} b_k^* q_k$ with $b_k^*=\lfloor (5/6) a_{k+1} \rfloor$. Noting that $F_k(N^*)=0$ for all $k$, from Proposition \ref{PNUNprop} we get
\[ \log P_N (\alpha) - \log P_{N^*}(\alpha) = \log U_N - \log U_{N^*} - \sum_{k=k_0}^{K-1} F_k(N) + O_T \left( \sum_{k=1}^K \frac{1}{a_k} \right) + O_{\alpha} (1) , \]
where $F_k(N) \ge 0$ for all $k$, and $F_k(N) =0$ for all $k$ such that $b_k \le (1-\delta_T ) a_{k+1}$.

Let us now successively apply projections to $N$ in two rounds, as described in Section \ref{regoptsection}: the first round consists of regularizing projections in increasing order of the indices $k_0 \le k \le K-1$, and the second round consists of optimizing projections in increasing order of the indices. This way $N$ is transformed into the integer $\sum_{k=0}^{k_0-1} b_k q_k+\sum_{k=k_0}^{K-1} b_k^*q_k$. Since $U_N$ does not depend on the first $k_0$ Ostrowski coefficients, Proposition \ref{regoptprojections} allows us to write $\log U_{N^*} - \log U_N = \sum_{k=k_0}^{K-1} (d_k^{\mathrm{reg}}(N) + d_k^{\mathrm{opt}}(N))$, and hence
\[ \log P_N (\alpha) - \log P_{N^*}(\alpha)  = - \sum_{k=k_0}^{K-1} (d_k^{\mathrm{reg}}(N) + d_k^{\mathrm{opt}}(N) + F_k(N)) + O_T \left( \sum_{k=1}^K \frac{1}{a_k} \right) + O_{\alpha} (1) . \]
Here $d_k^{\mathrm{reg}}(N)$ resp.\ $d_k^{\mathrm{opt}}(N)$ describe the effect of the regularizing resp.\ optimizing projection with respect to the index $k$, and by Proposition \ref{regoptprojections} satisfy the following for all $k_0 \le k \le K-1$:
\begin{enumerate}
\item[(i)] if $b_k \le (1-\delta_T)a_{k+1}$, then $d_k^{\mathrm{reg}}(N) =0$;
\item[(ii)] if $b_k>(1-\delta_T ) a_{k+1}$, then $d_k^{\mathrm{reg}}(N) \ge -O_T (1/a_{k+1}) -O_{\alpha}(1/q_k)$;
\item[(iii)] if $b_k \le (1-\delta_T )a_{k+1}$, then
\[ \begin{split} d_k^{\mathrm{opt}} = &a_{k+1} \int_{b_k/a_{k+1}}^{b_k^*/a_{k+1}} \log |2 \sin (\pi x)| \, \mathrm{d}x \\ &+ O_T \left( \frac{|b_k-b_k^*|}{a_{k+1}} +I_{\{ b_k \le 0.01 a_{k+1} \}} \log a_{k+1} \right) +O_{\alpha}(1/q_k) ; \end{split} \]
\item[(iv)]
\[ d_k^{\mathrm{opt}} (N) \ge 0.2326 \frac{(b_k-b_k^*)^2}{a_{k+1}} - O_T (1/a_{k+1}) - O_{\alpha}(1/q_k) ; \]
\item[(v)] if $b_k=b_k^*$, then $d_k^{\mathrm{opt}}(N)=0$.
\end{enumerate}
Note that if $b_k>(1-\delta_T )a_{k+1}$, then $b_k$ is first projected to $b_k^{**}$ as defined in \eqref{bk**def}, and then to $b_k^*$. Proposition \ref{regoptprojections} (ii) thus yields property (iv) with $b_k^{**}$ in place of $b_k$. Choosing $\delta_T$ small enough, property (iv) also holds as stated with an arbitrarily smaller numerical constant. Observe also that the special value of $b_k^{**}$ in the case $a_{k+1}=2$ does not cause any problem.

It follows that we can introduce small error terms $\xi_k(N)=O_T (1/a_{k+1}) + O_{\alpha}(1/q_k)$ such that $d_k(N):= d_k^{\mathrm{reg}}(N) + d_k^{\mathrm{opt}}(N)+F_k(N)+\xi_k(N)$ satisfies $d_k (N) \ge 0.2326(b_k-b_k^*)^2 / a_{k+1}$ for all $k_0 \le k \le K-1$ with equality if $b_k=b_k^*$, and also
\[ d_k(N) = a_{k+1} \int_{b_k/a_{k+1}}^{b_k^*/a_{k+1}} \log |2 \sin (\pi x)| \, \mathrm{d}x + O_T \left( \frac{|b_k-b_k^*|}{a_{k+1}} +I_{\{ b_k \le 0.01 a_{k+1} \}} \log a_{k+1} \right) \]
for all $k_0 \le k \le K-1$ such that $b_k \le (1-\delta_T ) a_{k+1}$. Since the contribution of the error terms $\xi_k(N)$ is negligible, we also have
\[ \log P_N(\alpha ) - \log P_{N^*}(\alpha) = -\sum_{k=k_0}^{K-1} d_k(N) + O \left( \sum_{k=1}^K \frac{1}{a_k} \right) + O_{\alpha} (1) . \]
Finally, introduce $d_k(N)$, $0 \le k <k_0$ in any way which satisfy the desired properties, and observe that $\sum_{k=0}^{k_0-1} d_k(N)=O_{\alpha}(1)$. This concludes the proof of Theorem \ref{maxpntheorem}.
\end{proof}

\section{Proof of Theorem \ref{Lpnormtheorem}} \label{sec_proofs_2}

Note that \eqref{maxnorm} follows directly from Theorem \ref{maxpntheorem}, in particular from $d_k(N) \ge 0$; alternatively, it also follows from taking the limit in \eqref{Lpnorm} as $c \to \infty$. It will thus be enough to prove \eqref{Lpnorm}.

The main idea of the proof is that if we choose an integer $N$ randomly from the interval $0 \le N < q_K$, then its Ostrowski coefficients $b_k$, $0 \le k \le K-1$ are almost independent random variables, close to being uniformly distributed on $\{ 0,1,\dots, a_{k+1} \}$. As a Gaussian tail estimate will show, the coefficients $|b_k-b_k^*| \gg \sqrt{a_{k+1} \log a_{k+1}}$ have negligible contribution. By Theorem \ref{maxpntheorem}, and in particular \eqref{dkequ}, we thus have
\[ \begin{split} \sum_{N=0}^{q_K-1} P_N(\alpha )^c &\approx \sum_{b_0=0}^{a_1} \sum_{b_1=0}^{a_2} \cdots \sum_{b_{K-1}=0}^{a_K} P_{N^*}(\alpha )^c \exp \left( - \sum_{k=0}^{K-1} \frac{\pi \sqrt{3} c}{2} \cdot \frac{(b_k-b_k^*)^2}{a_{k+1}} \right) \\ &\approx P_{N^*} (\alpha )^c \prod_{k=0}^{K-1} \sum_{n \in \mathbb{Z}} \exp \left(- \pi \sqrt{3} c n^2 / (2a_{k+1}) \right) \\ &\approx P_{N^*} (\alpha )^c \prod_{k=0}^{K-1} \sqrt{\frac{2a_{k+1}}{\sqrt{3}c}} ,  \end{split} \]
which explains the main term in \eqref{Lpnorm}. We now give the formal proof.

\begin{proof}[Proof of Theorem \ref{Lpnormtheorem}] Let $c \ge 0.01$, and consider the intervals
\[ J_k = [0, 0.99 a_{k+1}] \cap \left[ b_k^*-r_k, b_k^*+ r_k \right] \cap \mathbb{Z} , \]
where
\[ r_k= 10 \sqrt{\frac{a_{k+1}}{c} \log \left( \frac{a_{k+1}}{c} +2 \right)} . \]

First, we prove the lower bound in \eqref{Lpnorm}. Note that for any $\underline{b}=(b_0,b_1, \dots, b_{K-1}) \in J_0 \times J_1 \times \cdots \times J_{K-1}$, the expression $N_{\underline{b}}=\sum_{k=0}^{K-1} b_k q_k$ is the Ostrowski expansion of an integer $0 \le N_{\underline{b}} < q_K$; moreover, we obtain each integer at most once. We wish to apply Theorem \ref{maxpntheorem} to $N_{\underline{b}}$, and simply discard all other integers in $[0,q_K)$ not of this form. Since for all $b_k \in J_k$ we have
\[ \frac{|b_k-b_k^*|}{a_{k+1}} + \frac{|b_k-b_k^*|^3}{a_{k+1}^2} \ll E_{k,c}:= \frac{\log^{1/2} (a_{k+1}/c+2)}{c^{1/2} a_{k+1}^{1/2}} +\frac{\log^{3/2} (a_{k+1}/c+2)}{c^{3/2} a_{k+1}^{1/2}} , \]
Theorem \ref{maxpntheorem} and \eqref{dkequ} yield
\[ P_{N_{\underline{b}}}(\alpha ) = P_{N^*}(\alpha ) \exp \left( - \sum_{k=0}^{K-1} \frac{\pi \sqrt{3}}{2} \cdot \frac{(b_k-b_k^*)^2}{a_{k+1}} \right) \exp \left( O_T \left( \sum_{k=1}^K \left( E_{k,c} + \frac{1}{a_k} \right) \right) + O_{\alpha}(1) \right) . \]
Consequently,
\[ \begin{split} \log \left( \sum_{N=0}^{q_K-1} P_N(\alpha )^c \right)^{1/c} \ge &\log \left( \sum_{\underline{b} \in J_0 \times J_1 \times \cdots \times J_{K-1}} P_{N_{\underline{b}}}(\alpha )^c \right)^{1/c} \\ \ge &\log P_{N^*} (\alpha) + \frac{1}{c} \log \sum_{\underline{b} \in J_0 \times J_1 \times \cdots \times J_{K-1}} \exp \left( - \sum_{k=0}^{K-1} \frac{\pi \sqrt{3} c}{2} \cdot \frac{(b_k-b_k^*)^2}{a_{k+1}} \right) \\ &- O_T \left( \sum_{k=1}^K \left( E_{k,c} + \frac{1}{a_k} \right) \right) - O_{\alpha} (1) .  \end{split} \]
The sum over $\underline{b}$ factors:
\[ \sum_{\underline{b} \in J_0 \times J_1 \times \cdots \times J_{K-1}} \exp \left( - \sum_{k=0}^{K-1} \frac{\pi \sqrt{3} c}{2} \cdot \frac{(b_k-b_k^*)^2}{a_{k+1}} \right) = \prod_{k=0}^{K-1} \sum_{b_k \in J_k} \exp \left( - \frac{\pi \sqrt{3}c}{2} \cdot \frac{(b_k-b_k^*)^2}{a_{k+1}} \right) . \]
Now let $A>0$ be a large universal constant, to be chosen. We establish a lower bound in the cases $a_{k+1}/c >A$ and $a_{k+1}/c \le A$ separately. First, assume that $a_{k+1}/c > A$. We then have
\[ r_k = 10 \sqrt{\frac{a_{k+1}}{c} \log \left( \frac{a_{k+1}}{c} +2 \right)} \le \frac{a_{k+1}}{10000 c} \le \frac{a_{k+1}}{100} , \]
provided that $A>0$ is large enough. In particular, $[b_k^*-r_k, b_k^* +r_k] \subset [0,0.99a_{k+1}]$. It is now easy to see that
\[ \begin{split} \sum_{b_k \in J_k} \exp \left( - \frac{\pi \sqrt{3}c}{2} \cdot \frac{(b_k-b_k^*)^2}{a_{k+1}} \right) &\ge \sum_{-r_k \le n \le r_k} \exp \left( -\pi \sqrt{3} c n^2 /(2a_{k+1}) \right) \\ &= \int_{-r_k}^{r_k} \exp \left( -\pi \sqrt{3} c x^2 /(2a_{k+1}) \right) \, \mathrm{d}x +O(1) \\ &= \sqrt{\frac{2a_{k+1}}{\sqrt{3}c}} +O(1) , \end{split} \]
and so
\[ \frac{1}{c} \log \sum_{b_k \in J_k} \exp \left( - \frac{\pi \sqrt{3}c}{2} \cdot \frac{(b_k-b_k^*)^2}{a_{k+1}} \right) \ge \frac{1}{2c} \log \frac{2a_{k+1}}{\sqrt{3}c} - O \left( \frac{1}{a_{k+1}} \right) \]
in the case $a_{k+1}/c >A$. If $a_{k+1}/c \le A$, then by noting that $b_k^* \in J_k$, the left hand side of the previous formula is non-negative, and thus the previous formula remains true. Altogether we obtain the lower bound
\[ \log \left( \sum_{N=0}^{q_K-1} P_N(\alpha )^c \right)^{1/c} \ge \log P_{N^*}(\alpha) + \frac{1}{2c} \sum_{k=1}^K \log \frac{2a_k}{\sqrt{3}c} - O_T \left( \sum_{k=1}^K \left( E_{k,c} +\frac{1}{a_k} \right) \right) - O_{\alpha} (1) . \]

Next, we prove the upper bound in \eqref{Lpnorm}. Applying Theorem \ref{maxpntheorem} to all $0 \le N <q_K$, we get
\[ \log \left( \sum_{N=0}^{q_K-1} P_N(\alpha )^c \right)^{1/c} = \log P_{N^*} + \frac{1}{c} \log \sum_{N=0}^{q_K-1} \exp \left( -c \sum_{k=0}^{K-1} d_k(N) \right) + O_T \left( \sum_{k=1}^K \frac{1}{a_k} \right) + O_{\alpha}(1) . \]
For the sake of readability, note that $0.2326 >1/5$. Letting
\[ g_k(x) = \left\{ \begin{array}{ll} \frac{\pi \sqrt{3}}{2} \cdot \frac{(x-b_k^*)^2}{a_{k+1}} & \textrm{if } x \in J_k, \\ \frac{(x-b_k^*)^2}{5a_{k+1}} & \textrm{if } x \not\in J_k, \end{array} \right. \]
we have $d_k(N) \ge g_k(b_k) - O_T (E_{k,c})$. Hence by extending the range of summation,
\[ \begin{split} \log \left( \sum_{N=0}^{q_K-1} P_N(\alpha )^c \right)^{1/c} \le &\log P_{N^*} + \frac{1}{c} \log \sum_{\underline{b} \in \mathbb{Z}^K} \exp \left( -c \sum_{k=0}^{K-1} g_k(b_k) \right) \\ &+ O_T \left( \sum_{k=1}^K \left( E_{k,c} + \frac{1}{a_k} \right) \right) + O_{\alpha}(1) . \end{split} \]
The sum over $\underline{b}$ factors again:
\[ \frac{1}{c} \log \sum_{\underline{b}\in \mathbb{Z}^K} \exp \left( -c \sum_{k=0}^{K-1} g_k(b_k) \right) = \frac{1}{c} \sum_{k=0}^{K-1} \log \sum_{b_k \in \mathbb{Z}} \exp \left( -c g_k(b_k) \right) . \]
Now let $B>0$ be a large universal constant, to be chosen. Assume first, that $a_{k+1}/c>B$. Then, as before, $[b_k^*-r_k,b_k^*+r_k] \subset [0,0.99 a_{k+1}]$ provided that $B>0$ is large enough. Therefore
\[ \begin{split} & \sum_{b_k \in \mathbb{Z}} \exp (-cg_k (b_k) ) \\  \le & \sum_{|b_k-b_k^*| \le r_k} \exp \left( -\frac{\pi \sqrt{3}c}{2} \cdot \frac{(b_k-b_k^*)^2}{a_{k+1}} \right) + \sum_{|b_k-b_k^*| > r_k} \exp \left( -c \frac{(b_k-b_k^*)^2}{5a_{k+1}} \right) \\ \le & \int_{-\infty}^{\infty} \exp (-\pi \sqrt{3} c x^2/(2a_{k+1})) \, \mathrm{d}x + \int_{(-\infty, -r_k) \cup (r_k,\infty )} \exp (-c x^2/(5a_{k+1})) \, \mathrm{d}x +O(1) \\ = & \sqrt{\frac{2a_{k+1}}{\sqrt{3}c}} +O(1) ,  \end{split} \]
and consequently
\[  \frac{1}{c} \log \sum_{b_k \in \mathbb{Z}} \exp \left( -c g_k(b_k) \right) \le \frac{1}{2c} \log \frac{2a_{k+1}}{\sqrt{3}c} + O \left( \frac{1}{a_{k+1}} \right) . \]
If $a_{k+1}/c \le B$, then simply using $g_k(b_k) \ge (b_k-b_k^*)^2/(5a_{k+1})$ we similarly deduce that the left hand side of the previous formula is $\le O(1/c)$; consequently, the previous formula remains true. Altogether we obtain the upper bound
\[ \log \left( \sum_{N=0}^{q_K-1} P_N(\alpha )^c \right)^{1/c} \le \log P_{N^*}(\alpha) + \frac{1}{2c} \sum_{k=1}^K \log \frac{2a_k}{\sqrt{3}c} + O_T \left( \sum_{k=1}^K \left( E_{k,c} +\frac{1}{a_k} \right) \right) + O_{\alpha} (1) . \]
This concludes the proof of \eqref{Lpnorm}.
\end{proof}

\section{Proof of Theorem \ref{pn*theorem}} \label{sec_proofs_3}

In this section we estimate $P_{N^*}(\alpha)$. By Proposition \ref{PNUNprop} it is enough to consider $U_{N^*}$; in particular, we will need to estimate
\[ \begin{split} \log u_k (N^*) = &\sum_{b=1}^{b_k^*-1} \log |2 \sin (\pi (bq_k \| q_k \alpha \| + \varepsilon_k (N^*)))| + \sum_{b=0}^{b_k^*-1} V_k (b q_k \| q_k \alpha \| + \varepsilon_k (N^*)) \\ &+ \log (2 \pi (b_k^* q_k \| q_k \alpha \| + \varepsilon_k (N^*))) . \end{split} \]
The first sum can be handled with a straightforward application of a second order Euler--Maclaurin formula. Since it provides an elementary explanation for the appearance of the constant $\V$, we include a detailed proof in Lemma \ref{sumlogflemma} below. An estimate for the second sum follows from Lemma \ref{Vkxlemma} (iii); Theorem \ref{pn*theorem} will then be an immediate corollary.

We now give a formal proof. We will need the value of the integrals
\begin{equation}\label{improper1}
\int_1^{\infty} \frac{B_2(\{ x \})}{(x-5/6)^2} \, \mathrm{d}x = \frac{1}{3} - \log \frac{\Gamma (1/6)}{2^{5/6} 3^{1/3} \pi^{1/2}}
\end{equation}
and
\begin{equation}\label{improper2}
\int_1^{\infty} \frac{B_2(\{ x \})}{x^2} \, \mathrm{d}x = -\frac{11}{12} + \log (2^{1/2} \pi^{1/2}) ,
\end{equation}
where $B_2(x)=x^2/2-x/2+1/12$ is the second Bernoulli polynomial, and $\Gamma$ is the Gamma function. Indeed, by Stirling's formula for the Gamma function,
\[ \begin{split} \sum_{k=1}^n \log (k-5/6) = \log \frac{\Gamma (n+1/6)}{\Gamma (1/6)} = &(n-5/6) \log (n-5/6) -(n-5/6) \\ &+ \log \sqrt{2 \pi (n-5/6)} - \log \Gamma (1/6) + o(1) . \end{split} \]
On the other hand, applying a second order Euler--Maclaurin formula we get
\[ \begin{split} \sum_{k=1}^n \log (k-5/6) = &\int_1^n \log (x-5/6) \, \mathrm{d}x + \frac{\log (n-5/6) +\log (1/6)}{2} \\ &+ \frac{1/(n-5/6) - 6}{12} + \int_1^n \frac{B_2(\{ x \} )}{(x-5/6)^2} \, \mathrm{d}x , \end{split} \]
and by comparing the asymptotics as $n \to \infty$ in the previous formulas, \eqref{improper1} follows; the proof of \eqref{improper2} is analogous.
\begin{lem}\label{sumlogflemma} Let $N^*=\sum_{k=0}^{K-1} b_k^* q_k$ with $b_k^*=\lfloor (5/6)a_{k+1} \rfloor$. For any $0 \le k \le K-2$,
\[ \begin{split} \sum_{b=1}^{b_k^*-1} \log |2 \sin (\pi (b q_k \| q_k \alpha \| + \varepsilon_k (N^*)))| = &\frac{\V}{4 \pi q_k \| q_k \alpha \|} -\frac{1}{3} \log a_{k+1} - \log \frac{\Gamma (1/6)}{(2 \pi )^{5/6}} \\ &+ O \left( \frac{1}{a_{k+1}} + \frac{1+\log a_{k+1}}{a_{k+2}} \right) , \end{split} \]
whereas for $k=K-1$,
\[ \sum_{b=1}^{b_k^*-1} \log |2 \sin (\pi (b q_k \| q_k \alpha \| + \varepsilon_k (N^*)))| = \frac{\V}{4 \pi q_{K-1} \| q_{K-1} \alpha \|} +\frac{1}{2} \log a_K  + O \left( \frac{1}{a_K} \right) . \]
\end{lem}

\begin{proof} For the sake of readability, let $f(x) = |2 \sin (\pi x)|$ and $\varepsilon_k=\varepsilon_k (N^*)$. By the definition \eqref{varepsilonkdef} of $\varepsilon_k$ and the construction of $b_k^*$, for all $0 \le k \le K-2$ we have
\[ \begin{split} \varepsilon_k &= q_k \sum_{\ell =k+1}^{K-1} (-1)^{k+\ell} (5/6) a_{\ell +1} \| q_{\ell} \alpha \| + O \left( q_k \sum_{\ell =k+1}^{K-1} \| q_{\ell} \alpha \| \right) \\ &= -(5/6) q_k \| q_k \alpha \| +O(q_k \| q_{k+1} \alpha \|) , \end{split} \]
and in particular
\begin{equation}\label{varepsilonk5/6}
\frac{\varepsilon_k}{q_k \| q_k \alpha \|} = - \frac{5}{6} + O \left( \frac{1}{a_{k+2}} \right) ,
\end{equation}
whereas $\varepsilon_{K-1}=0$. Since $b_{k+1}^* \le (5/6)a_{k+2}$, Lemma \ref{epsilonklemma} also gives $q_k \| q_k \alpha \| + \varepsilon_k \ge 1/(18 a_{k+1})$.

Consider the function with first and second derivatives
\[ \begin{split} F(x) &= \log f(xq_k \| q_k \alpha \| + \varepsilon_k ), \\ F'(x) &= \pi \cot (\pi (x q_k \| q_k \alpha \| + \varepsilon_k ) ) q_k \| q_k \alpha \| , \\ F''(x) &= - \frac{\pi^2 q_k^2 \| q_k \alpha \|^2}{\sin^2 (\pi (x q_k \| q_k \alpha \| + \varepsilon_k ))} . \end{split} \]
Applying a second order Euler--Maclaurin formula, we get
\begin{equation}\label{eulermaclaurin}
\begin{split} \sum_{b=1}^{b_k^*-1} F(b) = &\int_1^{b_k^*-1} F(x) \, \mathrm{d}x + \frac{F(b_k^*-1)+F(1)}{2} + \frac{F'(b_k^*-1) - F'(1)}{12} \\ &- \int_1^{b_k^*-1} F''(x) B_2(\{ x \} ) \mathrm{d}x . \end{split}
\end{equation}
First, we estimate the main term
\[ \int_1^{b_k^*-1} F(x) \, \mathrm{d}x = \frac{1}{q_k \| q_k \alpha \|} \int_{q_k \| q_k \alpha \| + \varepsilon_k}^{(b_k^*-1) q_k \| q_k \alpha \| + \varepsilon_k} \log f (y) \, \mathrm{d}y . \]
Here by construction $(b_k^*-1) q_k \| q_k \alpha \| + \varepsilon_k = 5/6 + O(1/a_{k+1})$. Since $\log f(5/6)=0$, the error of replacing the upper limit of integration by $5/6$ is negligible:
\[ \frac{1}{q_k \| q_k \alpha \|} \left| \int_{(b_k^*-1) q_k \| q_k \alpha \| + \varepsilon_k}^{5/6} \log f(y) \mathrm{d}y \right| \ll \frac{1}{a_{k+1}} . \]
The effect of replacing the lower limit of integration by $0$ is
\[ \begin{split} \frac{1}{q_k \| q_k \alpha \|} \int_0^{q_k \| q_k \alpha \| + \varepsilon_k} \log f(y) \, \mathrm{d} y = &\frac{1}{q_k \| q_k \alpha \|} \int_0^{q_k \| q_k \alpha \| + \varepsilon_k} \log (2 \pi y) \, \mathrm{d}y + O \left( \frac{1}{a_{k+1}^2} \right) \\ = &\left( 1+\frac{\varepsilon_k}{q_k \| q_k \alpha \|} \right) \left( \log ( 2 \pi (q_k \| q_k \alpha \| +\varepsilon_k )) -1 \right) \\ &+ O \left( \frac{1}{a_{k+1}^2} \right) . \end{split} \]
From the previous three formulas and \eqref{varepsilonk5/6} it follows that the main term in \eqref{eulermaclaurin} is, for all $0 \le k \le K-2$,
\begin{equation}\label{eulerfirstterm1}
\int_1^{b_k^*-1} F(x) \, \mathrm{d}x = \frac{1}{q_k \| q_k \alpha \|} \int_0^{5/6} \log f(y) \, \mathrm{d}y + \frac{1}{6}+\frac{1}{6} \log \frac{3 a_{k+1}}{\pi} + O \left( \frac{1}{a_{k+1}} + \frac{1+\log a_{k+1}}{a_{k+2}} \right) ,
\end{equation}
whereas for $k=K-1$,
\begin{equation}\label{eulerfirstterm2}
\int_1^{b_{K-1}^*-1} F(x) \, \mathrm{d}x = \frac{1}{q_{K-1} \| q_{K-1} \alpha \|} \int_0^{5/6} \log f(y) \, \mathrm{d}y +1+ \log \frac{a_K}{2 \pi} + O \left( \frac{1}{a_K} \right) .
\end{equation}
Next, consider the second and third terms in \eqref{eulermaclaurin}. It is easy to see that $F(b_k^*-1) \ll 1/a_{k+1}$ and $F'(b_k^*-1) \ll 1/a_{k+1}$. Further,
\[ F(1)= \log f(q_k \| q_k \alpha \| + \varepsilon_k ) = \log (2 \pi (q_k \| q_k \alpha \| + \varepsilon_k )) + O \left( \frac{1}{a_{k+1}^2} \right) , \]
and using $\pi \cot (\pi y) =1/y +O(|y|)$ on $(0,5/6]$,
\[ F'(1) = \pi \cot (\pi (q_k \| q_k \alpha \| + \varepsilon_k )) q_k \| q_k \alpha \| = \frac{q_k \| q_k \alpha \|}{q_k \| q_k \alpha \| + \varepsilon_k} + O \left( \frac{1}{a_{k+1}^2} \right) . \]
Hence by \eqref{varepsilonk5/6}, for $0 \le k \le K-2$ we have
\begin{equation}\label{eulersecondterm1}
\frac{F(b_k^*-1)+F(1)}{2} + \frac{F'(b_k^*-1) - F'(1)}{12} =- \frac{1}{2} - \frac{1}{2} \log \frac{3a_{k+1}}{\pi} + O \left( \frac{1}{a_{k+1}} + \frac{1}{a_{k+2}} \right) ,
\end{equation}
whereas for $k=K-1$,
\begin{equation}\label{eulersecondterm2}
\frac{F(b_{K-1}^*-1)+F(1)}{2} + \frac{F'(b_{K-1}^*-1) - F'(1)}{12} = - \frac{1}{12} - \frac{1}{2} \log \frac{a_K}{2 \pi} + O \left( \frac{1}{a_K} \right) .
\end{equation}
Finally, consider the last term in \eqref{eulermaclaurin}. Using $\pi^2/\sin^2 (\pi y) = 1/y^2 +O(1)$ on $(0,5/6]$, 
\[ - \int_1^{b_k^*-1} F''(x) B_2 (\{ x \} ) \, \mathrm{d}x = \int_1^{b_k^*-1} \frac{B_2 (\{ x \} )}{\left( x+ \frac{\varepsilon_k}{q_k \| q_k \alpha \|} \right)^2} \, \mathrm{d}x +O \left( \frac{1}{a_{k+1}} \right) . \]
Therefore by \eqref{varepsilonk5/6} and the improper integrals \eqref{improper1} and \eqref{improper2}, for all $0 \le k \le K-2$ we have
\begin{equation}\label{eulerlastterm1}
- \int_1^{b_k^*-1} F''(x) B_2 (\{ x \} ) \, \mathrm{d}x = \frac{1}{3} - \log \frac{\Gamma (1/6)}{2^{5/6}3^{1/3}\pi^{1/2}} + O \left( \frac{1}{a_{k+1}} + \frac{1}{a_{k+2}} \right) ,
\end{equation}
whereas for $k=K-1$,
\begin{equation}\label{eulerlastterm2}
- \int_1^{b_{K-1}^*-1} F''(x) B_2 (\{ x \} ) \, \mathrm{d}x = -\frac{11}{12} + \log (2^{1/2} \pi^{1/2})  + O \left( \frac{1}{a_K} \right) .
\end{equation}

We have thus estimated all terms in \eqref{eulermaclaurin}. The claim for $0 \le k \le K-2$ follows from \eqref{eulerfirstterm1}, \eqref{eulersecondterm1} and \eqref{eulerlastterm1}, whereas the claim for $k=K-1$ follows from \eqref{eulerfirstterm2}, \eqref{eulersecondterm2} and \eqref{eulerlastterm2}.
\end{proof}

\begin{proof}[Proof of Theorem \ref{pn*theorem}] Let $\varepsilon_k = \varepsilon_k (N^*)$. Applying Proposition \ref{PNUNprop} and noting that $F_k(N^*)=0$ for all $k$ we get
\[ \begin{split} \log P_{N^*}(\alpha ) = &\sum_{k=0}^{K-1} \sum_{b=1}^{b_k^*-1} \log |2 \sin (\pi (bq_k \| q_k \alpha \| + \varepsilon_k ))| + \sum_{k=0}^{K-1} \sum_{b=0}^{b_k^*-1} V_k (b q_k \| q_k \alpha \| + \varepsilon_k ) \\ &+ \sum_{k=0}^{K-1} I_{\{ b_k^* \ge 1 \}} \log (2 \pi (b_k^* q_k \| q_k \alpha \| + \varepsilon_k )) + O_T \left( \sum_{k=1}^K \frac{1}{a_k} \right) + O_{\alpha} (1) . \end{split} \]
Note that $I_{\{ b_k^* \ge 1 \}} = I_{\{ a_{k+1} \ge 2 \}}$. The first sum was evaluated in Lemma \ref{sumlogflemma}. We can estimate the second sum by interpreting it as a Riemann sum and using Lemma \ref{Vkxlemma} (iii). Note that the endpoints are $\varepsilon_k =O(1/a_{k+1})$ and $(b_k^*-1)q_k \| q_k \alpha \| + \varepsilon_k =5/6+O(1/a_{k+1})$. Since the points $b q_k \| q_k \alpha \| + \varepsilon_k$, $b=0,1,\dots, b_k^*-1$ lie in the interval $[-(1-\frac{1}{e^T+2} ), 5/6]$, and since by Lemma \ref{Vkxlemma} the function $V_k$ is monotonically decreasing and satisfies $|V_k(x)|\ll_T (1+\log a_k)/a_{k+1}$ on this interval, we have
\[ \begin{split} \sum_{b=0}^{b_k^*-1} V_k (b q_k \| q_k \alpha \| + \varepsilon_k ) = &\frac{1}{q_k \| q_k \alpha \|} \int_0^{5/6} V_k (x) \, \mathrm{d}x + O_T \left( \frac{1+\log a_k}{a_{k+1}} \right) \\ = &\int_0^{5/6} \left( \log \frac{a_k}{2\pi} - \frac{\Gamma' (1+x)}{\Gamma (1+x)} \right) \, \mathrm{d} x \\ &+ O_T \left( \frac{1 + \log (a_{k-1}a_k)}{a_k} + \frac{1+\log a_k}{a_{k+1}} \right) + O_{\alpha} \left( \frac{1}{q_{k+1}} \right) \\ = &\frac{5}{6} \log \frac{a_k}{2 \pi} - \log \Gamma \left( 1+\frac{5}{6} \right) \\ &+ O_T \left( \frac{1 + \log (a_{k-1}a_k)}{a_k} + \frac{1+\log a_k}{a_{k+1}} \right) + O_{\alpha} \left( \frac{1}{q_{k+1}} \right) . \end{split} \]
Clearly,
\[ I_{\{ b_k^* \ge 1 \}} \log (2 \pi (b_k^* q_k \| q_k \alpha \| + \varepsilon_k )) = \log \left( \frac{5}{6} 2 \pi \right) + O_T \left( \frac{1}{a_{k+1}} \right) . \]
Summing over $0 \le k \le K-1$, from the previous two formulas and Lemma \ref{sumlogflemma} we get
\[ \begin{split} \log P_{N^*} (\alpha ) = &\frac{\V}{4 \pi} \sum_{k=0}^{K-1} \frac{1}{q_k \| q_k \alpha \|} -\frac{1}{3} \sum_{k=0}^{K-2} \log a_{k+1} + \frac{1}{2} \log a_K -K \log \frac{\Gamma (1/6)}{(2 \pi )^{5/6}} \\ &+ \frac{5}{6} \sum_{k=0}^{K-1} \log a_k + K \left( - \frac{5}{6} \log (2 \pi ) - \log \Gamma \left( 1+\frac{5}{6} \right) + \log \left( \frac{5}{6} 2\pi \right) \right) \\ &+O_T \left( \sum_{k=1}^{K-1} \frac{1 + \log (a_k a_{k+1})}{a_{k+1}} \right) + O_{\alpha}(1) . \end{split} \]
Observe that with a remarkable cancellation the coefficient of $K$ vanishes. Indeed, by Euler's reflection formula $\Gamma (x) \Gamma (1-x) = \pi /\sin (\pi x)$ we have $\Gamma (1/6) \Gamma (5/6) =2 \pi$, and hence
\begin{equation} \label{remark_canc}
- \log \frac{\Gamma (1/6)}{(2 \pi )^{5/6}} - \frac{5}{6} \log (2 \pi ) - \log \Gamma \left( 1 + \frac{5}{6} \right) + \log \left( \frac{5}{6} 2 \pi \right) =0. 
\end{equation}
The previous formula for $\log P_{N^*}(\alpha)$ thus simplifies to
\[ \log P_{N^*}(\alpha ) = \frac{\V}{4 \pi} \sum_{k=0}^{K-1} \frac{1}{q_k \| q_k \alpha \|} + \frac{1}{2} \sum_{k=1}^{K} \log a_k +O_T \left( \sum_{k=1}^{K-1} \frac{1 + \log (a_k a_{k+1})}{a_{k+1}} \right) + O_{\alpha}(1) . \]
Using e.g.\ property (iii) of continued fractions in Section \ref{sec_cf}, we see that here $1/(q_k \| q_k \alpha \|) = a_{k+1}+O(1/a_k+1/a_{k+2})$, and the claim follows.
\end{proof}

\section{Proof of Theorems \ref{quadraticlimitfunction} and \ref{wellapproximablelimitfunction}} \label{sec_proofs_4}

\subsection{Quadratic irrationals}\label{quadraticsubsection}

In this section we estimate the limit functions of $P_{q_k}(\alpha, (-1)^k x/q_k)$ for a given quadratic irrational $\alpha$. In order to make our estimates uniform on the interval of interest $(-1,1)$, we isolate the singularities at $x=-1$ and $1$. To this end, let us introduce the modified cotangent sum
\[ V_k^*(x):= \sum_{\substack{1 \leq n \leq q_k-1, \\ n \neq q_{k-1}, q_k-q_{k-1}}} \sin (\pi n \| q_k \alpha \| /q_k) \cot \left( \pi \frac{n(-1)^k p_k+x}{q_k} \right) . \]
Observe that by excluding $n=q_{k-1}$ resp.\ $n=q_k-q_{k-1}$, we avoid $n(-1)^k p_k \equiv -1 \pmod{q_k}$ resp.\ $n(-1)^k p_k \equiv 1 \pmod{q_k}$. In particular, $V_k^* (x)$ does not have a singularity on $(-2,2)$. The evaluation in Lemma \ref{Vkxlemma} (iii) has a perfect analogue for $V_k^*(x)$.
\begin{lem}\label{Vk*xlemma} Assume \eqref{logak/ak+1}. For any $k \ge 4$ and any $x \in (-2,2)$,
\[ \frac{V_k^*(x)}{q_k \| q_k \alpha \|} = \log \frac{a_k}{2 \pi} - \frac{\Gamma' (2+x)}{\Gamma (2+x)} +O \left( \frac{T+\log (a_{k-1}a_k)}{(2-|x|)a_k} \right) + O_{\alpha}(1/q_k) . \]
\end{lem}

\begin{proof} Following the proof of Lemma \ref{Vkxlemma} with obvious modifications, we get
\[ V_k^*(x) = \pi \| q_k \alpha \| \sum_{\substack{1 \leq n \leq q_k-1, \\ n \neq q_{k-1}, q_k-q_{k-1}}} \frac{n}{q_k} \cot \left( \pi \frac{n(-1)^k p_k+x}{q_k} \right) +O\left( \frac{\| q_k \alpha \|^3 q_k \log q_k}{2-|x|} \right) . \]
It remains to prove
\begin{align}
C_k^* (x) & := \sum_{\substack{1 \leq n \leq q_k-1, \\ n \neq q_{k-1}, q_k-q_{k-1}}} \frac{n}{q_k} \cot \left( \pi \frac{n p_k+(-1)^k x}{q_k} \right) \nonumber\\ & = \frac{(-1)^k q_k}{\pi} \left( \log \frac{a_k}{2\pi} -\frac{\Gamma' (2+x)}{\Gamma (2+x)} + O \left( \frac{T+\log (a_{k-1}a_k)}{(2-|x|)a_k} \right) \right) +O_{\alpha}(1). \label{ck*x}
\end{align}
From Lemma \ref{cotangentevaluation} we obtain
\[ \begin{split} C_k^*(0) &= \sum_{n=1}^{q_k-1} \frac{n}{q_k} \cot \left( \pi \frac{np_k}{q_k} \right) - \frac{q_{k-1}}{q_k} \cot \left( \pi \frac{(-1)^{k+1}}{q_k} \right) - \frac{q_k-q_{k-1}}{q_k} \cot \left( \pi \frac{(-1)^k}{q_k} \right) \\ &= \frac{(-1)^k q_k}{\pi} \left( \log \frac{a_k}{2 \pi} + \gamma -1+O \left( \frac{T+\log (a_{k-1}a_k)}{a_k} \right) \right) + O_{\alpha} (1) , \end{split} \]
where $\gamma =-\Gamma' (1)/\Gamma (1)$ is the Euler--Mascheroni constant. Note that we used $\cot (\pi /q_k) = q_k/\pi + O(1)$ and that $q_{k-1}/q_k \ll 1/a_k$ is negligible. Following the proof of Lemma \ref{cotangentevaluation} with obvious modifications (note that excluding $n=q_{k-1}$ and $n=q_k-q_{k-1}$ corresponds to excluding $a=\pm 1$), we get that the derivative of $C_k^*(x)$ satisfies
\[ {C_k^*}'(x) = \frac{(-1)^{k+1}q_k}{\pi} \sum_{a=2}^{\infty} \frac{1}{(a+x)^2} +O \left( \frac{q_k (1+\log a_k)}{(2-|x|)^2 a_k} \right) . \]
By integrating,
\[ \begin{split} C_k^*(x) - C_k^*(0) &= \frac{(-1)^{k+1}q_k}{\pi} \sum_{a=2}^{\infty} \left( \frac{1}{a} - \frac{1}{a+x} \right) + O \left( \frac{q_k (1+\log a_k)}{(2-|x|) a_k} \right) \\ &= \frac{(-1)^{k+1} q_k}{\pi} \left( \gamma -1+ \frac{\Gamma' (2+x)}{\Gamma (2+x)} +O \left( \frac{1+\log a_k}{(2-|x|) a_k} \right) \right) , \end{split} \]
and \eqref{ck*x} for general $x \in (-2,2)$ follows. 
\end{proof}

Next, let us introduce the appropriately modified version of $B_{k,M}(x)$ from \eqref{BNx}: for any integers $k \ge 1$ and $0 \le M <q_k$, let
\[ \begin{split} B_{k,M}^*(x):= \log \frac{P_M^* (\alpha , (-1)^k x/q_k)}{P_M^* (p_k/q_k , (-1)^kx/q_k)} - \sum_{\substack{1 \leq n \leq M, \\ n \neq q_{k-1}, q_k-q_{k-1}}} \sin (\pi \| q_k \alpha \| /q_k) \cot \left( \pi \frac{n(-1)^k p_k+x}{q_k} \right) , \end{split} \]
where
\[ P_M^*(\alpha , (-1)^k x/q_k) := \prod_{\substack{1 \leq n \leq M, \\ n \neq q_{k-1}, q_k-q_{k-1}}} |2 \sin (\pi (n \alpha + (-1)^k x/q_k))| , \]
and $P_M^*(p_k/q_k, (-1)^k x/q_k)$ is defined analogously.
\begin{prop}\label{modifiedtransferprop}
Let $k \ge 1$ and $0 \le M < q_k$ be integers, and assume that $q_k \| q_k \alpha \| \le 2(1-c_k)$ and $-2<x \le 2-\frac{q_k \| q_k \alpha \|}{1-c_k}$ with some $100/q_k^2 \le c_k<1$. Then
\[ -C \frac{\log (4/c_k)}{(2-|x|)^2 a_{k+1}^2} \le B_{k,M}^*(x) \le C \frac{1}{a_{k+1}^2 q_k} \]
with a universal constant $C>0$.
\end{prop}

\begin{proof} This is an obvious modification of the proof of Proposition \ref{transferprop} (i).
\end{proof}

\begin{proof}[Proof of Theorem \ref{quadraticlimitfunction}] Let $\alpha=[a_0;a_1, \dots, a_{k_0}, \overline{a_{k_0+1}, \dots, a_{k_0+p}}]$ be a quadratic irrational, and assume that $\max_{1 \le r \le p} (\log a_{k_0+r})/a_{k_0+r+1} \le T$ with some constant $T \ge 1$. From Corollary \ref{Bqk-1corollary} we deduce
\begin{equation}\label{pqkquadratic}
\begin{split} \log P_{q_k}(\alpha, (-1)^k x/q_k) = &\log \left( |2 \sin (\pi (\| q_k \alpha \| + x/q_k))| \frac{|\sin (\pi x)|}{|\sin (\pi x/q_k)|} \right) \\ &+ \sum_{n \in \{ q_{k-1}, q_k-q_{k-1} \}} \log \frac{|\sin (\pi (n \alpha + (-1)^k x/q_k))|}{|\sin (\pi (n p_k/q_k + (-1)^k x/q_k))|} \\ &+V_k^*(x)+B_{k,q_k-1}^*(x) . \end{split}
\end{equation}
Recall from the proof of Lemma \ref{epsilonklemma} that $q_k \| q_k \alpha \| \le 1-\frac{1}{e^T+2}$. Applying Proposition \ref{modifiedtransferprop} with $c_k=\frac{1}{e^T+2} \le \frac{1}{2}$ we thus obtain that for all $|x| \le \max \{ 1 , 2-2/a_{k+1} \}$ and all large enough $k$ (in terms of $T$),
\[ |B_{k,q_k-1}^*(x)| \ll \frac{T}{(2-|x|)^2 a_{k+1}^2} . \]
Applying Lemma \ref{Vk*xlemma}, formula \eqref{pqkquadratic} thus simplifies to
\begin{equation}\label{pqkquadratic2}
\begin{split} \log P_{q_k}(\alpha, (-1)^k x/q_k) = &\log \left( |2 \sin (\pi (\| q_k \alpha \| + x/q_k))| \frac{|\sin (\pi x)|}{|\sin (\pi x/q_k)|} \right) \\ &+ \sum_{n \in \{ q_{k-1}, q_k-q_{k-1} \}} \log \frac{|\sin (\pi (n \alpha + (-1)^k x/q_k))|}{|\sin (\pi (n p_k/q_k + (-1)^k x/q_k))|} \\ &+q_k \| q_k \alpha \| \left( \log \frac{a_k}{2 \pi} - \frac{\Gamma' (2+x)}{\Gamma (2+x)} \right) \\ &+ O \left( \frac{T+\log (a_{k-1}a_k)}{(2-|x|)a_k a_{k+1}} + \frac{T}{(2-|x|)^2 a_{k+1}^2} \right) + O_{\alpha}(1/q_k) . \end{split}
\end{equation}
We now let $k \to \infty$ along the arithmetic progression $p\mathbb{N} + k_0 + r$, and claim that every term in \eqref{pqkquadratic2} (except the first error term) converges. Indeed, we clearly have $q_k \| q_k \alpha \| \to C_r$ and $q_{k-1} \| q_k \alpha \| \to D_r$ with some constants $C_r, D_r >0$ depending on $\alpha$. The limit of the first term in \eqref{pqkquadratic2} is
\[ \log \left( |2 \sin (\pi (\| q_k \alpha \| + x/q_k))| \frac{|\sin (\pi x)|}{|\sin (\pi x/q_k)|} \right) \to \log \left( |2 \sin (\pi x)| \cdot \left| 1+\frac{C_r}{x} \right| \right) . \]
Using trigonometric identities, we once again write
\[ \frac{|\sin (\pi (n \alpha + (-1)^k x/q_k))|}{|\sin (\pi (n p_k/q_k + (-1)^k x/q_k))|} = |1+x_n+y_n| \]
with
\[ x_n := \cos (\pi n (\alpha -p_k/q_k )) -1= \cos (\pi n \| q_k \alpha \| /q_k) -1 \]
and
\[ y_n := \sin (\pi n \| q_k \alpha \| /q_k ) \cot (\pi (n (-1)^k p_k/q_k + x/q_k)) . \]
For both $n=q_{k-1}$ and $n=q_k-q_{k-1}$ we have $x_n \to 0$. For $n=q_{k-1}$,
\[ y_{q_{k-1}} = \sin (\pi q_{k-1} \| q_k \alpha \| /q_k) \cot (\pi (x-1)/q_k) \to \frac{D_r}{x-1}, \]
whereas for $n=q_k-q_{k-1}$,
\[ y_{q_k-q_{k-1}} = \sin (\pi (q_k-q_{k-1})\| q_k \alpha \| /q_k) \cot (\pi (x+1)/q_k) \to \frac{C_r-D_r}{x+1} . \]
Note that we used $\cot (\pi y)=1/(\pi y)+O(1)$ as $y \to 0$. Therefore the limit of the second term in \eqref{pqkquadratic2} is
\[ \sum_{n \in \{ q_{k-1}, q_k-q_{k-1} \}} \log \frac{|\sin (\pi (n \alpha + (-1)^k x/q_k))|}{|\sin (\pi (n p_k/q_k + (-1)^k x/q_k))|} \to \log \left| 1+\frac{D_r}{x-1} \right| + \log \left| 1+\frac{C_r-D_r}{x+1} \right| . \]
From \eqref{pqkquadratic2} we thus get that for all $|x| \le \max \{ 1,2-2/a_{k_0+r+1} \}$,
\[ \begin{split} \log G_{\alpha, r}(x)= &\log \left( |2 \sin (\pi x)| \cdot \left| 1+\frac{C_r}{x} \right| \right) + \log \left| 1+\frac{D_r}{x-1} \right| + \log \left| 1+\frac{C_r-D_r}{x+1} \right| \\ &+ C_r \left( \log \frac{a_{k_0+r}}{2 \pi} -\frac{\Gamma' (2+x)}{\Gamma (2+x)} \right)  \\ & 
+ O \left( \frac{T+\log (a_{k_0+r-1}a_{k_0+r})}{(2-|x|)a_{k_0+r} a_{k_0+r+1}} + \frac{T}{(2-|x|)^2 a_{k_0+r+1}^2} \right) , \end{split} \]
as claimed.
\end{proof}

\subsection{Well approximable irrationals}

\begin{proof}[Proof of Theorem \ref{wellapproximablelimitfunction}] Let $\alpha$ be such that $\sup_{k \ge 1} a_k = \infty$. It will be enough to prove that
\begin{equation}\label{locallyuniform}
P_{q_{k_m}} (\alpha, (-1)^{k_m} x/q_{k_m}) e^{-V_{k_m}(0)} \to |2 \sin (\pi x)| \qquad \textrm{locally uniformly on } \mathbb{R}
\end{equation}
for any increasing sequence of positive integers $k_m$ such that $a_{k_m+1} \to \infty$ as $m \to \infty$; recall that $V_k(x)$ was defined in \eqref{Vkxdef}. Indeed, under the stronger assumption
\[ \frac{1+\log \max_{1 \le \ell \le k_m} a_{\ell}}{a_{k_m+1}} \to 0 \qquad \textrm{as } m \to \infty \]
we have $V_{k_m}(0) \to 0$ by Lemma \ref{Vkxlemma} (ii), therefore \eqref{locallyuniform} holds without the factor $e^{-V_{k_m}(0)}$, as claimed. If in addition $(1+\log a_k)/a_{k+1} \to 0$ (in particular, $a_{k+1} \to \infty$), then $V_k(0) \to 0$ by Lemma \ref{Vkxlemma} (iii), and \eqref{locallyuniform} follows without the factor $e^{-V_{k_m}(0)}$ along the full sequence $k_m=m$, as claimed.

Fix a large integer $A>0$, and let us prove that the convergence in \eqref{locallyuniform} is uniform on $[-A,A]$. Let
\[ S=S_k= \{ 1 \le n \le q_k-1 \, : \, np_k \equiv a \pmod{q_k} \textrm{ with some integer } 0<|a| \le A \} , \]
and let us introduce the modified cotangent sum
\[ V_k^{**}(x)= \sum_{\substack{1 \leq n \leq q_k-1, \\ n \not\in S}} \sin (\pi n \| q_k \alpha \| /q_k) \cot \left( \pi \frac{n(-1)^k p_k+x}{q_k} \right) . \]
Note that $V_k^{**}(x)$ does not have a singularity on $(-A-1,A+1)$. Following the steps in Section \ref{quadraticsubsection} with obvious modifications (cf.\ \eqref{pqkquadratic}), we get
\[ \begin{split} \log P_{q_k}(\alpha, (-1)^k x/q_k) = &\log \left( |2 \sin (\pi (\| q_k \alpha \| + x/q_k))| \frac{|\sin (\pi x)|}{|\sin (\pi x/q_k)|} \right) \\ &+ \sum_{n \in S} \log \frac{|\sin (\pi (n \alpha + (-1)^k x/q_k))|}{|\sin (\pi (n p_k/q_k + (-1)^k x/q_k))|} \\ &+V_k^{**}(x)+B_{k,q_k-1}^{**}(x) . \end{split} \]
Here $B_{k,q_k-1}^{**}(x)$ is the perfect analogue of $B_{k,q_k-1}^*(x)$ in \eqref{pqkquadratic}, and satisfies
\[ |B_{k,q_k-1}^{**}(x)| \ll \frac{1}{a_{k+1}^2}, \qquad x \in [-A,A], \]
by an obviously modified form of Proposition \ref{modifiedtransferprop} with $c_k=1/2$. Following the steps in the proof of Lemma \ref{Vkxlemma}, it is easy to see that the derivative of $V_k^{**}(x)$ satisfies $|{V_k^{**}}' (x)| \ll 1/a_{k+1}$ on $[-A,A]$. Therefore for any $x \in [-A,A]$,
\[ \begin{split} V_k^{**}(x) &= V_k^{**}(0) + O \left( \frac{A}{a_{k+1}} \right) \\ &= V_k(0) - \sum_{n \in S} \sin (\pi n \| q_k \alpha \| /q_k) \cot \left( \pi \frac{n(-1)^k p_k}{q_k} \right) + O \left( \frac{A}{a_{k+1}} \right) \\ &=V_k(0) + O \left( \sum_{0<|a| \le A} \| q_k \alpha \| \left| \cot \left( \pi \frac{a}{q_k} \right) \right| + \frac{A}{a_{k+1}} \right) \\ &= V_k(0) + O \left( \frac{A}{a_{k+1}} \right) . \end{split} \]
By the previous three formulas and the usual trigonometric identities,
\[ \begin{split} P_{q_k}(\alpha, (-1)^k x/q_k) e^{-V_k(0)} = |2 \sin (\pi x)| \frac{|\sin (\pi (\| q_k \alpha \| + x/q_k))|}{|\sin (\pi x/q_k)|} \left( \prod_{n \in S} |1+x_n+y_n| \right) e^{O \left( A/a_{k+1} \right)} \end{split} \]
uniformly on $[-A,A]$, where
\[ x_n := \cos (\pi n (\alpha -p_k/q_k )) -1= \cos (\pi n \| q_k \alpha \| /q_k) -1 \]
and
\[ y_n := \sin (\pi n \| q_k \alpha \| /q_k ) \cot (\pi (n (-1)^k p_k/q_k + x/q_k)) . \]
To see \eqref{locallyuniform}, it will thus be enough to prove that
\begin{equation}\label{wellapproxenough}
|2 \sin (\pi x)| \frac{|\sin (\pi (\| q_k \alpha \| + x/q_k))|}{|\sin (\pi x/q_k)|} \prod_{n \in S} |1+x_n+y_n| \to |2 \sin (\pi x)| \quad \textrm{uniformly on } [-A,A]
\end{equation}
along any subsequence $k=k_m$ such that $a_{k_m+1} \to \infty$.

First, let $x \in [-A,A] \backslash \bigcup_{a=-A}^A (a-1/100,a+1/100)$. Then
\[ \frac{|\sin (\pi (\| q_k \alpha \| + x/q_k))|}{|\sin (\pi x/q_k)|} \sim 1+\frac{q_k \| q_k \alpha \|}{x} , \]
as well as 
$$
|x_n| \ll \| q_k \alpha \|^2 \ll 1/a_{k+1} \qquad \text{and} \qquad |y_n| \ll_A 1/a_{k+1},
$$
all uniformly in $x$. Hence
\begin{align*}
\frac{|\sin (\pi (\| q_k \alpha \| + x/q_k))|}{|\sin (\pi x/q_k)|} \prod_{n \in S} |1+x_n+y_n|  & \sim \left( 1+O \left( \frac{1}{a_{k+1}} \right) \right) \prod_{0<|a| \le A} \left| 1+O_A \left( \frac{1}{a_{k+1}} \right) \right| \\
& =1+O_A \left( \frac{1}{a_{k+1}} \right), 
\end{align*}
uniformly in $x$. This shows that the convergence in \eqref{wellapproxenough} is indeed uniform on $[-A,A] \backslash \bigcup_{a=-A}^A (a-1/100,a+1/100)$. Next, let $x \in (a-1/100, a+1/100)$ with some $0<|a| \le A$. Then
\[ \begin{split} |2 \sin (\pi x)| \prod_{n \in S} |1+x_n+y_n| &= |2 \sin (\pi x)| \cdot \left| 1+O \left( \frac{1}{|x-a| a_{k+1}} \right) \right| \prod_{\substack{0<|a'| \le A \\ a' \neq a}} \left( 1+O_A \left( \frac{1}{a_{k+1}} \right) \right) \\ &= |2 \sin (\pi x)| + O_A \left( \frac{1}{a_{k+1}} \right), \end{split} \]
uniformly in $x \in (a-1/100, a+1/100)$; indeed, $|\sin (\pi x)|/|x-a| \ll 1$ follows from the fact that $|\sin (\pi x)|$ has a zero at every integer. Therefore the convergence in \eqref{wellapproxenough} is also uniform on $(a-1/100, a+1/100)$. Finally, let $x \in (-1/100, 1/100)$. Then
\[ |2 \sin (\pi x)| \frac{|\sin (\pi (\| q_k \alpha \| + x/q_k))|}{|\sin (\pi x/q_k)|} \sim |2 \sin (\pi x)| \left( 1+\frac{q_k \| q_k \alpha \|}{x} \right) = |2 \sin (\pi x)| + O \left( \frac{1}{a_{k+1}} \right), \]
uniformly in $x \in (-1/100, 1/100)$. Therefore the convergence in \eqref{wellapproxenough} is uniform on $(-1/100, 1/100)$. This finishes the proof of \eqref{wellapproxenough}.
\end{proof}

\section*{Acknowledgements}

CA is supported by the Austrian Science Fund (FWF), projects F-5512, I-3466, I-4945 and Y-901. BB is supported by FWF project Y-901.

\end{document}